\documentclass{article}

\newtheorem{theorem}{Theorem}[section]

\newtheorem{conjecture}[theorem]{Conjecture}
\newtheorem{lemma}[theorem]{Lemma}
\newtheorem{proposition}[theorem]{Proposition}

\evensidemargin\oddsidemargin
\begin{document}

\title{On Unimodality of Independence Polynomials of some Well-Covered Trees\thanks{%
A preliminary version of this paper was presented at Second Haifa Workshop
on Interdisciplinary Applications, of Graph Theory, Combinatorics and
Algorithms June 17-20, 2002.}}
\author{Vadim E. Levit and Eugen Mandrescu \\
Department of Computer Science\\
Holon Academic Institute of Technology\\
52 Golomb Str., P.O. Box 305\\
Holon 58102, ISRAEL}
\date{}
\maketitle

\begin{abstract}
The \textit{stability number }$\alpha (G)$ of the graph $G$ is the size of a
maximum stable set of $G$. If $s_{k}$ denotes the number of stable sets of
cardinality $k$ in graph $G$, then $I(G;x)=\sum\limits_{k=0}^{\alpha
(G)}s_{k}x^{k}$ is the \textit{independence polynomial} of $G$ (I. Gutman
and F. Harary, 1983). In 1990, Y. O. Hamidoune proved that for any \textit{%
claw-free graph} $G$ (a graph having no induced subgraph isomorphic to $%
K_{1,3}$), $I(G;x)$ is unimodal, i.e., there exists some $k\in
\{0,1,...,\alpha (G)\}$ such that 
\[
s_{0}\leq s_{1}\leq ...\leq s_{k}\geq s_{k+1}\geq ...\geq s_{\alpha (G)}. 
\]
Y. Alavi, P. J. Malde, A. J. Schwenk and P. Erd\"{o}s (1987) asked whether
for trees (or perhaps forests) the independence polynomial is unimodal. J.
I. Brown, K. Dilcher and R. J. Nowakowski (2000) conjectured that $I(G;x)$
is unimodal for any \textit{well-covered graph} $G$ (a graph whose all
maximal independent sets have the same size). V. E. Levit and E. Mandrescu
(1999) demonstrated that every well-covered tree can be obtained as a join
of a number of well-covered spiders, where a \textit{spider} is a tree
having at most one vertex of degree at least three.\textit{\ }

In this paper we show that the independence polynomial of any well-covered
spider is unimodal. In addition, we introduce some graph transformations
respecting independence polynomials. They allow us to reduce several types
of well-covered trees to claw-free graphs, and, consequently, to prove that
their independence polynomials are unimodal.\newline

\textbf{key words:}\textit{\ stable set, independence polynomial, unimodal
sequence, well-covered tree, claw-free graph.}
\end{abstract}

\section{Introduction}

Throughout this paper $G=(V,E)$ is a simple (i.e., a finite, undirected,
loopless and without multiple edges) graph with vertex set $V=V(G)$ and edge
set $E=E(G).$ If $X\subset V$, then $G[X]$ is the subgraph of $G$ spanned by 
$X$. By $G-W$ we mean the subgraph $G[V-W]$, if $W\subset V(G)$. We also
denote by $G-F$ the partial subgraph of $G$ obtained by deleting the edges
of $F$, for $F\subset E(G)$, and we write shortly $G-e$, whenever $F$ $%
=\{e\} $. The \textit{neighborhood} of a vertex $v\in V$ is the set $%
N_{G}(v)=\{w:w\in V$ \ \textit{and} $vw\in E\}$, and $N_{G}[v]=N_{G}(v)\cup
\{v\}$; if there is no ambiguity on $G$, we use $N(v)$ and $N[v]$,
respectively. If $N(v)$ induces a complete graph in $G$, then $v$ is a 
\textit{simplicial} vertex of $G$. A simplicial vertex is \textit{pendant}
if its neighborhood contains only one vertex, and an edge is \textit{pendant}
if at least one of its endpoints is a pendant vertex. $%
K_{n},P_{n},C_{n},K_{n_{1},n_{2},...,n_{p}}$ denote respectively, the
complete graph on $n\geq 1$ vertices, the chordless path on $n\geq 1$
vertices, the chordless cycle on $n\geq 3$ vertices, and the complete $p$%
-partite graph on $n_{1}+n_{2}+...+n_{p}$ vertices.

The \textit{disjoint union} of the graphs $G_{1},G_{2}$ is the graph $%
G=G_{1}\amalg G_{2}$ having as a vertex set the disjoint union of $%
V(G_{1}),V(G_{2})$, and as an edge set the disjoint union of $%
E(G_{1}),E(G_{2})$. In particular, $\amalg nG$ denotes the disjoint union of 
$n>1$ copies of the graph $G$. If $G_{1},G_{2}$ are disjoint graphs, then
their \textit{Zykov sum}, (Zykov, \cite{Zykov}, \cite{Zykov1}), is the graph 
$G_{1}+G_{2}$ with 
\begin{eqnarray*}
V(G_{1}+G_{2}) &=&V(G_{1})\cup V(G_{2}), \\
E(G_{1}+G_{2}) &=&E(G_{1})\cup E(G_{2})\cup \{v_{1}v_{2}:v_{1}\in
V(G_{1}),v_{2}\in V(G_{2})\}.
\end{eqnarray*}

As usual, a \textit{tree} is an acyclic connected graph. A tree having at
most one vertex of degree $\geq 3$ is called a \textit{spider}, \cite
{HedLaskar}, or an \textit{aster}, \cite{GorMcDonOrlYung}.

A \textit{stable} set in $G$ is a set of pairwise non-adjacent vertices. A
stable set of maximum size will be referred to as a \textit{maximum stable
set} of $G$, and the \textit{stability number }of $G$, denoted by $\alpha
(G) $, is the cardinality of a maximum stable set in $G$. Let $s_{k}$ be the
number of stable sets in $G$ of cardinality $k,k\in \{1,...,\alpha (G)\}$.
The polynomial 
\[
I(G;x)=\sum\limits_{k=0}^{\alpha (G)}s_{k}x^{k},s_{0}=1, 
\]
is called the \textit{independence polynomial} of $G$, (Gutman and Harary, 
\cite{GutHar}).

A number of general properties of the independence polynomial of a graph are
presented in \cite{GutHar} and \cite{Arocha}. As important examples, we
mention the following: 
\begin{eqnarray*}
I(G_{1}\amalg G_{2};x) &=&I(G_{1};x)\cdot I(G_{2};x), \\
I(G_{1}+G_{2};x) &=&I(G_{1};x)+I(G_{2};x)-1.
\end{eqnarray*}

A finite sequence of non-negative real numbers $%
\{a_{0},a_{1},a_{2},...,a_{n}\}$ is said to be \textit{unimodal} if there is
some $k\in \{0,1,...,n\}$, called the \textit{mode} of the sequence, such
that 
\[
0\leq a_{0}\leq a_{1}\leq ...\leq a_{k}\geq a_{k+1}\geq ...\geq a_{n}. 
\]
The mode is \textit{unique} if $a_{k-1}<a_{k}>a_{k+1}$.

Unimodal sequences occur in many areas of mathematics, including algebra,
combinatorics, and geometry (see Brenti, \cite{Brenti}, and Stanley, \cite
{Stanley}). In the context of our paper, for instance, if $a_{i}$ denotes
the number of ways to select a subset of $i$ independent edges (a matching
of size $i$) in a graph, then the sequence of these numbers is unimodal
(Schwenk, \cite{Schwenk}). As another example, if $a_{i}$ denotes the number
of dependent $i$-sets of a graph $G$ (sets of size $i$ that are not stable),
then the sequence of $\left\{ a_{i}\right\} _{i=0}^{n}$ is unimodal
(Horrocks, \cite{Horrocks}).

A polynomial $P(x)=a_{0}+a_{1}x+a_{2}x^{2}+...+a_{n}x^{n}$ is called \textit{%
unimodal} if the sequence of its coefficients is unimodal. For instance, the
independence polynomial of $K_{n}$ is unimodal, as $I(K_{n};x)=1+nx$.
However, the independence polynomial of the graph $G=K_{100}+\amalg 3K_{6}$
is not unimodal, since $I(G;x)=1+\mathbf{118}x+108x^{2}+\mathbf{206}x^{3}$
(for another examples, see Alavi et al \cite{AlMalSchErdos}). Moreover, in 
\cite{AlMalSchErdos} it is shown that for any permutation $\sigma $ of $%
\{1,2,...,\alpha \}$ there exists a graph $G$, with $\alpha (G)=\alpha $,
such that $s_{\sigma (1)}<s_{\sigma (2)}<...<s_{\sigma (\alpha )}$, i.e.,
there are graphs for which $s_{1},s_{2},...,s_{\alpha }$ is as ''shuffled''
as we like.

A graph $G$ is called \textit{well-covered} if all its maximal stable sets
have the same cardinality, (Plummer, \cite{Plum}). In particular, a tree $T$
is well-covered if and only if $T=K_{1}$ or it has a perfect matching
consisting of pendant edges (Ravindra, \cite{Ravindra}).

The roots of the independence polynomial of well-covered graphs are
investigated by Brown et al in \cite{Brown}. It is shown that for any
well-covered graph $G$ there is a well-covered graph $H$ with $\alpha
(G)=\alpha (H)$ such that $G$ is an induced subgraph of $H$, where all the
roots of $I(H;x)$ are simple and real. As it is also mentioned in \cite
{Brown}, a root of independence polynomial of a graph (not necessarily
well-covered) of smallest modulus is real, and there are well-covered graphs
whose independence polynomials have non-real roots. Moreover, it is easy to
check that the complete $n$-partite graph $G=K_{\alpha ,\alpha ,...,\alpha }$
is well-covered, $\alpha (G)=\alpha $, and its independence polynomial,
namely $I(G;x)=n(1+x)^{\alpha }-(n-1)$, has only one real root, whenever $%
\alpha $ is odd, and exactly two real roots, for any even $\alpha $. In
other words, the theorem of Newton (stating that if a polynomial with
positive coefficients has only real roots, then its coefficients form a
unimodal sequence) does not help in proving the following conjecture.

\begin{conjecture}
\cite{Brown} The independence polynomial of any well-covered graph is
unimodal.
\end{conjecture}

\begin{figure}[h]
\setlength{\unitlength}{1cm}%
\begin{picture}(5,2)\thicklines
  
  \multiput(3,0)(1,0){4}{\circle*{0.29}}
  \multiput(3,1)(1,0){4}{\circle*{0.29}}
  \multiput(3,2)(1,0){2}{\circle*{0.29}}
  \put(3,0){\line(1,0){1}} 
  \put(3,1){\line(1,0){3}}
  \put(3,2){\line(1,0){1}}
  \put(4,0){\line(0,1){2}} 
  \put(5,0){\line(0,1){1}}
  \put(6,0){\line(0,1){1}}

\put(2.3,1){\makebox(0,0){$T_{1}$}}

  \multiput(8,0)(1,0){5}{\circle*{0.29}}
  \multiput(8,1)(1,0){5}{\circle*{0.29}}
  \multiput(8,2)(1,0){2}{\circle*{0.29}}
  \put(8,0){\line(1,0){1}} 
  \put(8,1){\line(1,0){4}}
  \put(8,2){\line(1,0){1}}
  \put(9,0){\line(0,1){2}} 
  \put(10,0){\line(0,1){1}}
  \put(11,0){\line(0,1){1}}
  \put(12,0){\line(0,1){1}}  
  \put(7.3,1){\makebox(0,0){$T_{2}$}}
  
 \end{picture}
\caption{Two well-covered trees.}
\label{fig57}
\end{figure}
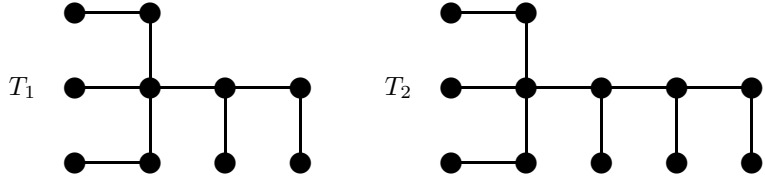

The claw-graph $K_{1,3}$ (see Figure \ref{fig7}) is a non-well-covered tree
and $I(K_{1,3};x)=1+4x+3x^{2}+x^{3}$ is unimodal, but has non-real roots.
The trees $T_{1},T_{2}$ in Figure \ref{fig57} are well-covered, and their
independence polynomials are respectively 
\begin{eqnarray*}
I(T_{1};x) &=&(1+x)^{2}\cdot (1+2x)\cdot (1+6x+7x^{2}) \\
&=&1+10x+36x^{2}+60x^{3}+47x^{4}+14x^{5}, \\
I(T_{2};x) &=&(1+6x+10x^{2}+5x^{3})^{2}-x^{2}(1+x)^{3}(1+2x) \\
&=&1+12x+55x^{2}+125x^{3}+151x^{4}+93x^{5}+23x^{6},
\end{eqnarray*}
which are both unimodal, while only for the first is true that all its roots
are real. Hence, Newton's theorem is not useful in verifying the following
conjecture, even for the particular case of well-covered\emph{\ }trees.

\begin{conjecture}
\cite{AlMalSchErdos} Independence polynomials of trees are unimodal.
\end{conjecture}

A graph is called \textit{claw-free} if it has no induced subgraph
isomorphic to $K_{1,3}$. There are non-claw-free graphs whose independence
polynomials are unimodal,\emph{\ }e.g., the $n$-star $K_{1,n},n\geq 3$. The
following result of Hamidoune will be used in the sequel.

\begin{theorem}
\cite{Hamidoune}\label{th2} The independence polynomial of a claw-free graph
is unimodal.
\end{theorem}

As a simple application of this statement, one can easily see that
independence polynomials of paths and cycles are unimodal. In \cite{Arocha},
Arocha shows that 
\[
I(P_{n};x)=F_{n+1}(x),\ and\ I(C_{n},x)=F_{n-1}(x)+2xF_{n-2}(x),
\]
where $F_{n}(x),n\geq 0$, are \textit{Fibonacci polynomials}, i.e., the
polynomials defined recursively by 
\[
F_{0}(x)=1,F_{1}(x)=1,F_{n}(x)=F_{n-1}(x)+xF_{n-2}(x).
\]
Based on this recurrence, one can deduce that 
\[
F_{n}(x)={n \choose 0}+{n-1 \choose 1}x + {n-2 \choose 2}x^{2}+...+
{{\lceil n/2\rceil } \choose {\lfloor n/2\rfloor }}x^{\lfloor n/2\rfloor },
\]
(for example, see Riordan, \cite{Riordan}, where this polynomial is
discussed as a special kind of rook polynomials). It is amusing that the
unimodality of the polynomial $F_{n}(x)$, which may be not so trivial to
establish directly, follows now immediately from Theorem \ref{th2}, since
any $P_{n}$ is claw-free. Let us notice that for $n\geq 5,P_{n}$ is not
well-covered. 
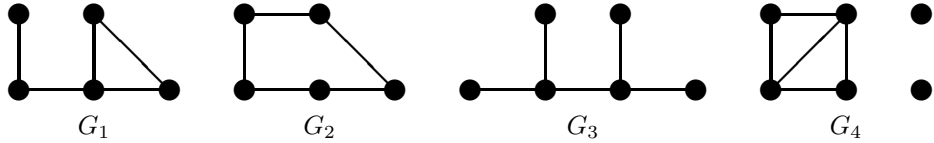
\begin{figure}[h]
\setlength{\unitlength}{1cm}%
\begin{picture}(5,1.7)\thicklines
  
  \multiput(0.5,0.5)(1,0){3}{\circle*{0.29}}
  \multiput(0.5,1.5)(1,0){2}{\circle*{0.29}}
  \put(0.5,0.5){\line(1,0){2}} 
  \put(0.5,0.5){\line(0,1){1}}
  \put(1.5,0.5){\line(0,1){1}}
  \put(1.5,1.5){\line(1,-1){1}}
  \put(1.5,0){\makebox(0,0){$G_{1}$}}

  \multiput(3.5,0.5)(1,0){3}{\circle*{0.29}}
  \multiput(3.5,1.5)(1,0){2}{\circle*{0.29}}
  \put(3.5,0.5){\line(1,0){2}} 
  \put(3.5,0.5){\line(0,1){1}}
  \put(3.5,1.5){\line(1,0){1}}
  \put(4.5,1.5){\line(1,-1){1}}
  \put(4.5,0){\makebox(0,0){$G_{2}$}}

  \multiput(6.5,0.5)(1,0){4}{\circle*{0.29}}
  \multiput(7.5,1.5)(1,0){2}{\circle*{0.29}}
  \put(6.5,0.5){\line(1,0){3}} 
  \put(7.5,0.5){\line(0,1){1}}
  \put(8.5,0.5){\line(0,1){1}}
  \put(8,0){\makebox(0,0){$G_{3}$}}

  \multiput(10.5,0.5)(1,0){3}{\circle*{0.29}}
  \multiput(10.5,1.5)(1,0){3}{\circle*{0.29}}
  \put(10.5,0.5){\line(1,0){1}} 
  \put(10.5,0.5){\line(0,1){1}}
  \put(10.5,0.5){\line(1,1){1}}
  \put(10.5,1.5){\line(1,0){1}}
  \put(11.5,0.5){\line(0,1){1}}

  \put(11.5,0){\makebox(0,0){$G_{4}$}}
  
  \end{picture}
\caption{Two pairs of non-isomorphic graphs $G_{1},G_{2}$ and $G_{3},G_{4}$
satisfying $I(G_{1};x)=I(G_{2};x)$ and $I(G_{3};x)=I(G_{4};x)$.}
\label{fig63}
\end{figure}

Clearly, any two isomorphic graphs have the same independence polynomial.
The converse is not generally true. For instance, while $%
I(G_{1};x)=I(G_{2};x)=1+5x+5x^{2}$, the well-covered graphs $G_{1}$ and $%
G_{2}$ are non-isomorphic (see Figure \ref{fig63}).

In addition, the graphs $G_{3},G_{4}$ in Figure \ref{fig63}, have identical
independence polynomials $I(G_{3};x)=I(G_{4};x)=1+6x+10x^{2}+6x^{3}+x^{4}$,
while $G_{3}$ is a tree, and $G_{4}$ is not connected and has cycles.

However, if $I(G;x)=1+nx,n\geq 1$, then $G$ is isomorphic to $K_{n}$. Figure 
\ref{fig7} gives us a source of some more examples of such uniqueness.
Namely, the figure presents all the graphs of size four with the stability
number equal to three. A simple check shows that their independence
polynomials are different: 
\begin{eqnarray*}
I(K_{1,3};x) &=&1+4x+3x^{2}+x^{3}, \\
I(G_{1};x) &=&1+4x+5x^{2}+2x^{3}, \\
I(G_{2};x) &=&1+4x+4x^{2}+x^{3}.
\end{eqnarray*}
In other words, if the independence polynomials of two graphs (one from
Figure \ref{fig7} and an arbitrary one) coincide, then these graphs are
exactly the same up to isomorphism.

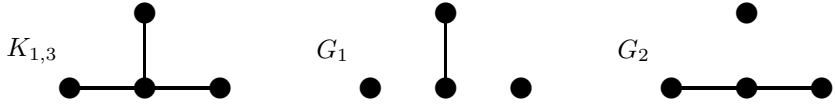
\begin{figure}[h]
\setlength{\unitlength}{1cm}%
\begin{picture}(5,1.2)\thicklines
  
  \multiput(2,0)(1,0){3}{\circle*{0.29}}
  \put(3,1){\circle*{0.29}}
  \put(2,0){\line(1,0){2}} 
  \put(3,0){\line(0,1){1}}
  \put(1.5,0.5){\makebox(0,0){$K_{1,3}$}}

  \multiput(6,0)(1,0){3}{\circle*{0.29}}
  \put(7,1){\circle*{0.29}}
  \put(7,0){\line(0,1){1}}
  \put(5.5,0.5){\makebox(0,0){$G_{1}$}}
  
  \multiput(10,0)(1,0){3}{\circle*{0.29}}
  \put(11,1){\circle*{0.29}}
  \put(10,0){\line(1,0){2}}
  \put(9.5,0.5){\makebox(0,0){$G_{2}$}}
  \end{picture}
\caption{$\alpha (K_{1,3})=\alpha (G_{1})=\alpha (G_{2})=3$.}
\label{fig7}
\end{figure}

Let us mention that the equality $I(G_{1};x)=I(G_{2};x)$ implies 
\[
\left| V(G_{1})\right| =s_{1}=\left| V(G_{2})\right| \ and\ %
\left| E(G_{1})\right| =\frac{s_{1}^{2}-s_{1}}{2}-s_{2}=\left|
E(G_{2})\right| . 
\]
Consequently, if $G_{1},G_{2}$ are connected, $I(G_{1};x)=I(G_{2};x)$ and
one of them is a tree, then the other must be a tree, as well.

In this paper we show that the independence polynomial of any well-covered
spider is unimodal. In addition, we introduce some graph transformations
respecting independence polynomials. They allow us to reduce several types
of well-covered trees to claw-free graphs, and, consequently, to prove that
their independence polynomials are unimodal.

\section{Preliminary results}

Let us notice that if the product of two polynomials is unimodal, this is
not a guaranty for the unimodality of at least one of the factors. For
instance, we have

\[
\begin{array}{l}
I(K_{100}+\amalg 3K_{6};x)\cdot I(K_{100}+\amalg 3K_{6};x)=(1+\mathbf{118}%
x+108x^{2}+\mathbf{206}x^{3})^{2} \\ 
=1+236x+14140x^{2}+25900x^{3}+\mathbf{60280}x^{4}+\allowbreak
44496x^{5}+42436x^{6}.
\end{array}
\]

The converse is also true: the product of two unimodal polynomials is not
necessarily unimodal. As an example, we see that:

\[
\begin{array}{l}
I(K_{100}+\amalg 3K_{7};x)\cdot I(K_{100}+\amalg 3K_{7};x)=(1+121x+147x^{2}+%
\mathbf{343}x^{3})^{2} \\ 
=\allowbreak 1+242x+14935x^{2}+36260x^{3}+\mathbf{104615}x^{4}+\allowbreak
100842x^{5}+\mathbf{117649}x^{6}.
\end{array}
\]

However, if one of them is of degree one, we show that their product is
still unimodal.

\begin{lemma}
\label{lem3}If $R_{n}$ is a unimodal polynomial, then $R_{n}\cdot R_{1}$ is
unimodal for any polynomial $R_{1}$.
\end{lemma}

\setlength {\parindent}{0.0cm}\textbf{Proof.} Let $%
R_{n}(x)=a_{0}+a_{1}x+a_{2}x^{2}+...+a_{n}x^{n}$ be a unimodal polynomial
and $R_{1}(x)=b_{0}+b_{1}x$. \setlength
{\parindent}{3.45ex}Suppose that $a_{0}\leq a_{1}\leq ...\leq a_{k}\geq
a_{k+1}\geq ...\geq a_{n}$ and $b_{0}\leq b_{1}$. Then, $R_{n}(x)\cdot
R_{1}(x)=a_{0}b_{0}+\sum\limits_{i=1}^{n}(a_{i}b_{0}+a_{i-1}b_{1})\cdot
x^{i}+a_{n}b_{1}\cdot x^{n+1}=\sum\limits_{i=0}^{n+1}c_{i}\cdot x^{i}$ and
we show that $R_{n}\cdot R_{1}$ is unimodal, with the mode $m$, where 
\begin{equation}
c_{m}=\max \{c_{k},c_{k+1}\}=\max
\{a_{k}b_{0}+a_{k-1}b_{1},a_{k+1}b_{0}+a_{k}b_{1}\}.  \label{ModLe}
\end{equation}

Firstly, $a_{0}b_{0}\leq a_{1}b_{0}+a_{0}b_{1}$ because $a_{0}b_{0}\leq \max
\{a_{1}b_{0},a_{0}b_{1}\}$ and $0\leq \min \{a_{1}b_{0},a_{0}b_{1}\}$.
Secondly, $a_{i-1}\leq a_{i}\leq a_{i+1}$ are true for any $i\in
\{1,...,k-1\}$, and these assure that $a_{i}b_{0}+a_{i-1}b_{1}\leq
a_{i+1}b_{0}+a_{i}b_{1}$. Further, $a_{i-1}\geq a_{i}\geq a_{i+1}$ are valid
for any $i\in \{k+1,...,n-1\}$, which imply that $a_{i}b_{0}+a_{i-1}b_{1}%
\geq a_{i+1}b_{0}+a_{i}b_{1}$. Finally, $a_{n}b_{0}+a_{n-1}b_{1}\geq
a_{n}b_{1}$, since $a_{n-1}\geq a_{n}$.

Similarly, we can show that $R_{n}\cdot R_{1}$ is unimodal, whenever $%
b_{0}>b_{1}$. \rule{2mm}{2mm}\newline

The following proposition constitutes an useful tool in computing
independence polynomials of graphs and also in finding recursive formulae
for independence polynomials of various classes of graphs.

\begin{proposition}
\cite{GutHar}, \cite{HoedeLi}\label{prop1} Let $G=(V,E)$ be a graph, $w\in
V,uv\in E$ and $U\subset V$ be such that $G[U]$ is a complete subgraph of $G$%
. Then the following equalities hold:

(i) $I(G;x)=I(G-w;x)+x\cdot I(G-N[w];x)$;

(ii) $I(G;x)=I(G-U;x)+x\cdot \sum\limits_{v\in U}I(G-N[v];x)$;

(iii) $I(G;x)=I(G-uv;x)-x^{2}\cdot I(G-N(u)\cup N(v);x)$.
\end{proposition}

The \textit{edge-join} of two disjoint graphs $G_{1},G_{2}$ is the graph $%
G_{1}\ominus G_{2}$ obtained by adding an edge joining two vertices
belonging to $G_{1},G_{2}$, respectively. If the two vertices are $v_{i}\in
V(G_{i}),i=1,2$, then by $(G_{1};v_{1})\ominus (G_{2};v_{2})$ we mean
the graph $G_{1}\ominus G_{2}$.

\begin{lemma}
\label{lem1}Let $G_{i}=(V_{i},E_{i}),i=1,2$, be two well-covered graphs and $%
v_{i}\in V_{i},i=1,2$, be simplicial vertices in $G_{1},G_{2}$,
respectively, such that $N_{G_{i}}[v_{i}],i=1,2$, contains at least another
simplicial vertex. Then the following assertions are true:

(i) $G=(G_{1};v_{1})$ $\ominus (G_{2};v_{2})$ is well-covered and $%
\alpha (G)=\alpha (G_{1})+\alpha (G_{2})$;

(ii) $I(G;x)=I(G_{1};x)\cdot I(G_{2};x)-x^{2}\cdot
I(G_{1}-N_{G_{1}}[v_{1}];x)\cdot I(G_{2}-N_{G_{2}}[v_{2}];x)$.
\end{lemma}

\setlength {\parindent}{0.0cm}\textbf{Proof.} (i) Let $S_{1},S_{2}$ be
maximum stable sets in $G_{1},G_{2}$, respectively. Since $G_{1},G_{2}$ are
well-covered, we may assume that $v_{i}\notin S_{i},i=1,2$. Hence, $%
S_{1}\cup S_{2}$ is stable in $G$ and any maximum stable set $A$ of $G$ has $%
\left| A\cap V_{1}\right| \leq \left| S_{1}\right| ,\left| A\cap
V_{2}\right| \leq \left| S_{2}\right| $, and consequently we obtain: 
\[
\left| S_{1}\right| +\left| S_{2}\right| =\left| S_{1}\cup S_{2}\right| \leq
\left| A\right| =\left| A\cap V_{1}\right| +\left| A\cap V_{2}\right| \leq
\left| S_{1}\right| +\left| S_{2}\right| , 
\]
i.e., $\alpha (G_{1})+\alpha (G_{2})=\alpha (G)$.%
\setlength
{\parindent}{3.45ex}

Let $B$ be a stable set in $G$ and $B_{i}=B\cap V_{i},i=1,2$. Clearly, at
most one of $v_{1},v_{2}$ may belong to $B$. Since $G_{1},G_{2}$ are
well-covered, there exist $S_{1},S_{2}$ maximum stable sets in $G_{1},G_{2}$%
, respectively, such that $B_{1}\subseteq S_{1},B_{2}\subseteq S_{2}$.

\textit{Case 1.} $v_{1}\in B$ (similarly, if $v_{2}\in B$), i.e., $v_{1}\in
B_{1}$. If $v_{2}\notin S_{2}$, then $S_{1}\cup S_{2}$ is a maximum stable
set in $G$ such that $B\subset S_{1}\cup S_{2}$. Otherwise, let $w$ be the
other simplicial vertex belonging to $N_{G_{2}}[v_{2}]$. Then $%
S_{3}=S_{2}\cup \{w\}-\{v_{2}\}$ is a maximum stable set in $G_{2}$, that
includes $B_{2}$, because $B_{2}\subseteq S_{2}-\{v_{2}\}$. Hence, $%
S_{1}\cup S_{3}$ is a maximum stable set in $G$ such that $B\subset
S_{1}\cup S_{3}$.

\textit{Case 2.} $v_{1},v_{2}\notin B$. If $v_{1},v_{2}\in S_{1}\cup S_{2}$,
then as above, $S_{1}\cup (S_{2}\cup \{w\}-\{v_{2}\})$ is a maximum stable
set in $G$ that includes $B$. Otherwise, $S_{1}\cup S_{2}$ is a maximum
stable set in $G$ such that $B\subset S_{1}\cup S_{2}$.

Consequently, $G=(G_{1};v_{1})$ $\ominus \ (G_{2};v_{2})$ is well-covered.

(ii) Using Proposition \ref{prop1}(iii), we obtain that
\begin{eqnarray*}
I(G;x) &=&I(G-v_{1}v_{2};x)-x^{2}\cdot I(G-N_{G}(v_{1})\cup N_{G}(v_{2});x)
\\
&=&I(G_{1};x)\cdot I(G_{2};x)-x^{2}\cdot I(G_{1}-N_{G_{1}}(v_{1});x)\cdot
I(G_{2}-N_{G_{2}}(v_{2});x),
\end{eqnarray*}
which completes the proof. \rule{2mm}{2mm}\newline

\begin{figure}[h]
\setlength{\unitlength}{1cm}%
\begin{picture}(5,2)\thicklines
  
  \multiput(2.5,0.5)(1,0){6}{\circle*{0.29}}
  \multiput(3,1.5)(2,0){3}{\circle*{0.29}}
  \put(9.5,0.5){\circle*{0.29}}
  \put(10.5,0.5){\circle*{0.29}}
  \put(10,1.5){\circle*{0.29}}
  \put(2.5,0.5){\line(1,0){5}} 
  \put(2.5,0.5){\line(1,2){0.5}}
  \put(4.5,0.5){\line(1,2){0.5}}
  \put(6.5,0.5){\line(1,2){0.5}} 
  \put(9.5,0.5){\line(1,2){0.5}}
  \put(3.5,0.5){\line(-1,2){0.5}}
  \put(5.5,0.5){\line(-1,2){0.5}}
  \put(7.5,0.5){\line(-1,2){0.5}} 
  \put(10.5,0.5){\line(-1,2){0.5}}
  \put(9.5,0.5){\line(1,0){1}}
  \multiput(7.5,0.5)(0.125,0){16}{\circle*{0.07}}
  \put(2.4,0){\makebox(0,0){$v_{1}$}}
  \put(3.5,0){\makebox(0,0){$v_{2}$}}
  \put(4.5,0){\makebox(0,0){$v_{4}$}}
  \put(5.5,0){\makebox(0,0){$v_{5}$}}
  \put(6.5,0){\makebox(0,0){$v_{7}$}}
  \put(7.5,0){\makebox(0,0){$v_{8}$}}
  \put(9.5,0){\makebox(0,0){$v_{3n-2}$}}
  \put(10.8,0){\makebox(0,0){$v_{3n-1}$}}
  \put(3,1.9){\makebox(0,0){$v_{3}$}}
  \put(5,1.9){\makebox(0,0){$v_{6}$}}
  \put(7,1.9){\makebox(0,0){$v_{9}$}}
  \put(10,1.9){\makebox(0,0){$v_{3n}$}}
  
  \end{picture}
\caption{The graph $\bigtriangleup _{n}=K_{3}\ominus (n-1)K_{3}$.}
\label{fig2}
\end{figure}
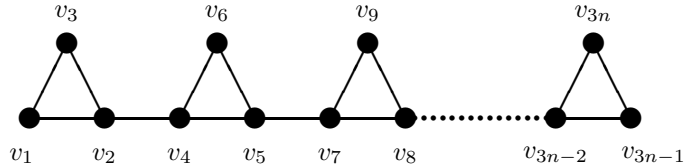

By $\bigtriangleup _{n}$ we denote the graph $\ominus nK_{3}$ defined
as $\bigtriangleup _{n}=K_{3}\ominus (n-1)K_{3},n\geq 1$, (see Figure 
\ref{fig2}). $\bigtriangleup _{0}$ denotes the empty graph.

\begin{proposition}
\label{prop4}The following assertions are true:

(i) for any $n\geq 1$, the graphs $\bigtriangleup _{n},K_{2}\ominus
\bigtriangleup _{n}$ are well-covered;

(ii) $I(\bigtriangleup _{n};x)$ is unimodal for any $n\geq 1$, and 
\[
I(\bigtriangleup _{n};x)=(1+3x)\cdot I(\bigtriangleup _{n-1};x)-x^{2}\cdot
I(\bigtriangleup _{n-2};x),n\geq 2, 
\]

where $I(\bigtriangleup _{0};x)=1,I(\bigtriangleup _{1};x)=1+3x$;

(iii) $I(K_{2}\ominus \bigtriangleup _{n};x)$ is unimodal for any $%
n\geq 1$, and 
\[
I(K_{2}\ominus \bigtriangleup _{n};x)=(1+2x)\cdot I(\bigtriangleup
_{n};x)-x^{2}\cdot I(\bigtriangleup _{n-1};x). 
\]
\end{proposition}

\setlength {\parindent}{0.0cm}\textbf{Proof.} (i) We show, by induction on $%
n $, that $\bigtriangleup _{n}$ is well-covered.%
\setlength
{\parindent}{3.45ex} Clearly, $\bigtriangleup _{1}=K_{3}$ is well-covered.
For $n\geq 2$ we have $\bigtriangleup _{n}=(\bigtriangleup
_{1};v_{2})\ominus (\bigtriangleup _{n-1};v_{4})$, (see Figure \ref
{fig2}). Hence, according to Lemma \ref{lem1}, $\bigtriangleup _{n}$ is
well-covered, because $v_{2},v_{3}$ and $v_{4},v_{6}$ are simplicial
vertices in $\bigtriangleup _{1},\bigtriangleup _{n-1}$, respectively.

Therefore, $\bigtriangleup _{n}$ is well-covered for any $n\geq 1$.

(ii) If $e=v_{2}v_{4}$ and $n\geq 2$, then according to Proposition \ref
{prop1}(iii), we obtain that 
\begin{eqnarray*}
I(\bigtriangleup _{n};x) &=&I(\bigtriangleup _{n}-e;x)-x^{2}\cdot
I(\bigtriangleup _{n}-N(v_{2})\cup N(v_{4});x) \\
&=&I(K_{3};x)\cdot I(\bigtriangleup _{n-1};x)-x^{2}\cdot I(\bigtriangleup
_{n-2};x) \\
&=&(1+3x)\cdot I(\bigtriangleup _{n-1};x)-x^{2}\cdot I(\bigtriangleup
_{n-2};x).
\end{eqnarray*}
In addition, $I(\bigtriangleup _{n};x)$ is unimodal by Theorem \ref{th2},
because $\bigtriangleup _{n}$ is claw-free.

\begin{figure}[h]
\setlength{\unitlength}{1cm}%
\begin{picture}(5,2)\thicklines

  \multiput(2,0.5)(1,0){7}{\circle*{0.29}}
  \multiput(3.5,1.5)(2,0){3}{\circle*{0.29}}
  \put(2,1.5){\circle*{0.29}}
  \put(10,0.5){\circle*{0.29}}
  \put(11,0.5){\circle*{0.29}}
  \put(10.5,1.5){\circle*{0.29}}
  \put(2,0.5){\line(1,0){6}}
  \put(2,0.5){\line(0,1){1}}
  \put(3,0.5){\line(1,2){0.5}}
  \put(5,0.5){\line(1,2){0.5}}
  \put(7,0.5){\line(1,2){0.5}} 
  \put(10,0.5){\line(1,2){0.5}}
  \put(4,0.5){\line(-1,2){0.5}}
  \put(6,0.5){\line(-1,2){0.5}}
  \put(8,0.5){\line(-1,2){0.5}} 
  \put(11,0.5){\line(-1,2){0.5}}
  \put(10,0.5){\line(1,0){1}}
  \multiput(8,0.5)(0.125,0){16}{\circle*{0.07}}
  \put(3.85,1.6){\makebox(0,0){$v_{3}$}}
  \put(4,0){\makebox(0,0){$v_{2}$}}
  \put(3,0){\makebox(0,0){$v_{1}$}}
  \put(1.6,0.5){\makebox(0,0){$u_{2}$}}
  \put(2.4,1.5){\makebox(0,0){$u_{1}$}}
  \put(5,0){\makebox(0,0){$v_{4}$}}
  \put(6,0){\makebox(0,0){$v_{5}$}}
  \put(7,0){\makebox(0,0){$v_{7}$}}
  \put(8,0){\makebox(0,0){$v_{8}$}}
  \put(10,0){\makebox(0,0){$v_{3n-2}$}}
  \put(11.7,0.5){\makebox(0,0){$v_{3n-1}$}}
  \put(5.85,1.6){\makebox(0,0){$v_{6}$}}
  \put(7.85,1.6){\makebox(0,0){$v_{9}$}}
  \put(10.95,1.6){\makebox(0,0){$v_{3n}$}}

  \end{picture}
\caption{The graph $K_{2}\ominus \bigtriangleup _{n}$.}
\label{fig32}
\end{figure}
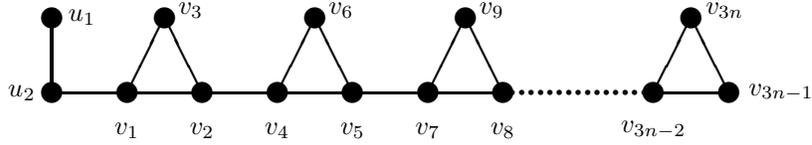

(iii) Let us notice that both $K_{2}$ and $\bigtriangleup _{n}$ are
well-covered. The graph $K_{2}\ominus \bigtriangleup
_{n}=(K_{2};u_{2})\ominus (\bigtriangleup _{n};v_{1})$ is well-covered
according to Lemma \ref{lem1}, and $I(K_{2}\ominus \bigtriangleup
_{n};x)$ is unimodal for any $n\geq 1$, by Theorem \ref{th2}, since $%
K_{2}\ominus \bigtriangleup _{n}$ is claw-free (see Figure \ref{fig32}).

In addition, applying Proposition \ref{prop1}(iii), we infer that

\begin{eqnarray*}
I(K_{2}\ominus \bigtriangleup _{n};x) &=&I(K_{2}\ominus
\bigtriangleup _{n}-u_{2}v_{1};x)-x^{2}\cdot I(K_{2}\ominus
\bigtriangleup _{n}-N(u_{2})\cup N(v_{1});x) \\
&=&I(K_{2};x)\cdot I(\bigtriangleup _{n};x)-x^{2}\cdot I(\bigtriangleup
_{n-1};x) \\
&=&(1+2x)\cdot I(\bigtriangleup _{n};x)-x^{2}\cdot I(\bigtriangleup
_{n-1};x),
\end{eqnarray*}
that completes the proof. \rule{2mm}{2mm}

\begin{lemma}
\label{lem2}Let $G_{i}=(V_{i},E_{i}),v_{i}\in V_{i},i=1,2$, and $%
P_{4}=(\{a,b,c,d\},\{ab,bc,cd\})$. Then the following assertions are true:

(i) $I(L_{1};x)=I(L_{2};x)$, where $L_{1}=(P_{4};b)\ominus (G_{1};v)$,
while $L_{2}$

\quad has $V(L_{2})=V(L_{1}),E(L_{2})=E(L_{1})\cup \{ac\}-\{cd\}$.

\quad If $G_{1}$ is claw-free and $v$ is simplicial in $G_{1}$, then $%
I(L_{1};x)$ is unimodal.

(ii) $I(G;x)=I(H;x)$, where $G=(G_{3};c)\ominus (G_{2};v_{2})$ and $%
G_{3}=(G_{1};v_{1})\ominus (P_{4};b)$,

\quad while $H$ has $V(H)=V(G),E(H)=E(G)\cup \{ac\}-\{cd\}$.

\quad If $G_{1},G_{2}$ are claw-free and $v_{1},v_{2}$ are simplicial in $%
G_{1},G_{2}$, respectively,

\quad then $I(G;x)$ is unimodal.
\end{lemma}

\begin{figure}[h]
\setlength{\unitlength}{1cm}%
\begin{picture}(5,2)\thicklines

  \multiput(2,0.5)(1,0){3}{\circle*{0.29}}
  \multiput(2,1.5)(1,0){2}{\circle*{0.29}}
  \put(2,0.5){\line(1,0){2}}
  \put(2,0.5){\line(0,1){1}}
  \put(3,0.5){\line(0,1){1}}
\multiput(3.7,0.1)(0.125,0){16}{\circle*{0.07}}
\multiput(3.7,1.9)(0.125,0){16}{\circle*{0.07}}
\multiput(3.7,0.1)(0,0.125){15}{\circle*{0.07}}
\multiput(5.7,0.1)(0,0.125){15}{\circle*{0.07}}
  \put(2,1.9){\makebox(0,0){$d$}}
  \put(2,0){\makebox(0,0){$c$}}
  \put(3,0){\makebox(0,0){$b$}}
  \put(3,1.9){\makebox(0,0){$a$}}
  \put(4.3,0.5){\makebox(0,0){$v$}}
  \put(4.7,1){\makebox(0,0){$G_{1}$}}

 \multiput(7.5,0.5)(1,0){3}{\circle*{0.29}}
  \multiput(7.5,1.5)(1,0){2}{\circle*{0.29}}
  \put(7.5,0.5){\line(1,0){2}}
  \put(7.5,0.5){\line(1,1){1}}
  \put(8.5,0.5){\line(0,1){1}}
  \multiput(9.2,1.9)(0.125,0){16}{\circle*{0.07}}
 \multiput(9.2,0.1)(0.125,0){16}{\circle*{0.07}}
 \multiput(9.2,0.1)(0,0.125){15}{\circle*{0.07}}
 \multiput(11.2,0.1)(0,0.125){15}{\circle*{0.07}}
 \put(8.5,1.9){\makebox(0,0){$a$}}
  \put(8.5,0){\makebox(0,0){$b$}}
  \put(7.5,0){\makebox(0,0){$c$}}
  \put(7.5,1.9){\makebox(0,0){$d$}}
  \put(10,0.5){\makebox(0,0){$v$}}
  \put(10.2,1){\makebox(0,0){$G_{1}$}}
 
  \end{picture}
\caption{The graphs $L_{1}=(P_{4};b)\ominus (G_{1};v)$ and $%
L_{2}=(K_{1}\amalg K_{3};b)\ominus (G_{1};v)$.}
\label{fig41}
\end{figure}
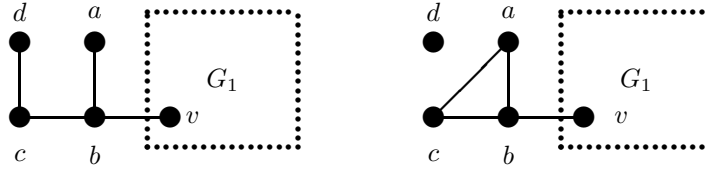

\setlength {\parindent}{0.0cm}\textbf{Proof.} (i) The graphs $%
L_{1}=(P_{4};b)\ominus (G_{1};v)$ and $L_{2}=(K_{1}\amalg
K_{3};b)\ominus (G_{1};v)$ are depicted in Figure \ref{fig41}.

Clearly, $I(P_{4};x)=I(K_{3}\amalg K_{1};x)=1+4x+3x^{2}$. 
\setlength
{\parindent}{3.45ex}By Proposition \ref{prop1}(iii), we obtain:

\begin{eqnarray*}
I(L_{1};x) &=&I(L_{1}-vb;x)-x^{2}\cdot I(L_{1}-N(v)\cup N(b);x) \\
&=&I(G_{1};x)\cdot I(P_{4};x)-x^{2}\cdot I(G_{1}-N_{G_{1}}\left[ v\right]
;x)\cdot I(\{d\};x).
\end{eqnarray*}

On the other hand, we get: 
\begin{eqnarray*}
I(L_{2};x) &=&I(L_{2}-vb;x)-x^{2}\cdot I(L_{2}-N(v)\cup N(b);x) \\
&=&I(G_{1};x)\cdot I(K_{3}\amalg K_{1};x)-x^{2}\cdot I(G_{1}-N_{G_{1}}\left[
v\right] ;x)\cdot I(\{d\};x).
\end{eqnarray*}

Consequently, the equality $I(L_{1};x)=I(L_{2};x)$ holds. If, in addition, $%
v $ is simplicial in $G_{1}$, and $G_{1}$ is claw-free, then $L_{2}$ is
claw-free, too. Theorem \ref{th2} implies that $I(L_{2};x)$ is unimodal,
and, hence, $I(L_{1};x)$ is unimodal, as well.

\begin{figure}[h]
\setlength{\unitlength}{1cm}%
\begin{picture}(5,2)\thicklines

  \multiput(1.5,0.5)(1,0){4}{\circle*{0.29}}
  \multiput(2.5,1.5)(1,0){2}{\circle*{0.29}}
  \put(1.5,0.5){\line(1,0){3}}
  \put(2.5,0.5){\line(0,1){1}}
  \put(3.5,0.5){\line(0,1){1}}
  \multiput(1.8,0.1)(-0.125,0){12}{\circle*{0.07}}
  \multiput(0.3,0.1)(0,0.125){15}{\circle*{0.07}}
  \multiput(1.8,1.9)(-0.125,0){12}{\circle*{0.07}}
  \multiput(1.8,0.1)(0,0.125){15}{\circle*{0.07}}
 \multiput(4.2,1.9)(0.125,0){12}{\circle*{0.07}}
 \multiput(4.2,0.1)(0.125,0){12}{\circle*{0.07}}
 \multiput(4.2,0.1)(0,0.125){15}{\circle*{0.07}}
 \multiput(5.7,0.1)(0,0.125){15}{\circle*{0.07}}
 \put(3.5,1.9){\makebox(0,0){$d$}}
  \put(3.5,0){\makebox(0,0){$c$}}
  \put(2.5,0){\makebox(0,0){$b$}}
  \put(2.5,1.9){\makebox(0,0){$a$}}
  \put(1.1,0.5){\makebox(0,0){$v_{1}$}}
  \put(0.8,1){\makebox(0,0){$G_{1}$}}
  \put(5,0.5){\makebox(0,0){$v_{2}$}}
  \put(5.2,1){\makebox(0,0){$G_{2}$}}  

 \multiput(8,0.5)(1,0){4}{\circle*{0.29}}
  \multiput(9,1.5)(1,0){2}{\circle*{0.29}}
  \put(8,0.5){\line(1,0){3}}
  \put(9,0.5){\line(0,1){1}}
  \put(10,0.5){\line(-1,1){1}}
  \multiput(8.3,0.1)(-0.125,0){12}{\circle*{0.07}}
  \multiput(6.8,0.1)(0,0.125){15}{\circle*{0.07}}
  \multiput(8.3,1.9)(-0.125,0){12}{\circle*{0.07}}
  \multiput(8.3,0.1)(0,0.125){15}{\circle*{0.07}}
 \multiput(10.7,1.9)(0.125,0){12}{\circle*{0.07}}
 \multiput(10.7,0.1)(0.125,0){12}{\circle*{0.07}}
 \multiput(10.7,0.1)(0,0.125){15}{\circle*{0.07}}
 \multiput(12.2,0.1)(0,0.125){15}{\circle*{0.07}}
 \put(10,1.9){\makebox(0,0){$d$}}
  \put(10,0){\makebox(0,0){$c$}}
  \put(9,0){\makebox(0,0){$b$}}
  \put(9,1.9){\makebox(0,0){$a$}}
  \put(7.6,0.5){\makebox(0,0){$v_{1}$}}
  \put(7.3,1){\makebox(0,0){$G_{1}$}}
  \put(11.5,0.5){\makebox(0,0){$v_{2}$}}
  \put(11.7,1){\makebox(0,0){$G_{2}$}}
 
  \end{picture}
\caption{The graphs $G$ and $H$ from Lemma 2.5(ii).}
\label{fig42}
\end{figure}
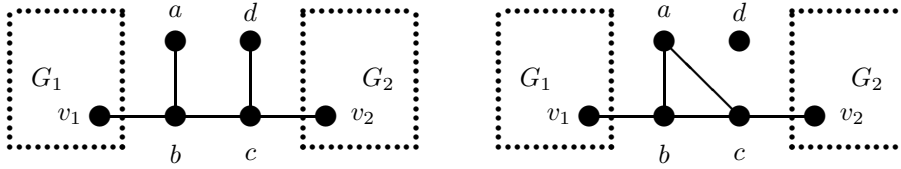

(ii) Figure \ref{fig42} shows the graphs $G$ and $H$.

According to Proposition \ref{prop1}(iii), we obtain: 
\begin{eqnarray*}
I(G;x) &=&I(G-v_{1}b;x)-x^{2}\cdot I(G-N(v_{1})\cup N(b);x) \\
&=&I(G_{1};x)\cdot I(G-G_{1};x)-x^{2}\cdot I(G_{1}-N_{G_{1}}\left[
v_{1}\right] ;x)\cdot I(G_{2};x)\cdot I(\left\{ d\right\} ;x).
\end{eqnarray*}

On the other hand, using again Proposition \ref{prop1}(iii), we get: 
\begin{eqnarray*}
I(H;x) &=&I(H-v_{1}b;x)-x^{2}\cdot I(H-N(v_{1})\cup N(b);x) \\
&=&I(G_{1};x)\cdot I(H-G_{1};x)-x^{2}\cdot I(G_{1}-N_{G_{1}}\left[
v_{1}\right] ;x)\cdot I(G_{2};x)\cdot I(\left\{ d\right\} ;x).
\end{eqnarray*}

Finally, let us observe that the equality $I(G-G_{1};x)=I(H-G_{1};x)$ holds
according to part (i).

Now, if $G_{1},G_{2}$ are claw-free and $v_{1},v_{2}$ are simplicial in $%
G_{1},G_{2}$, respectively, then $H$ is claw-free, and by Theorem \ref{th2},
its independence polynomial is unimodal. Consequently, $I(G;x)$ is also
unimodal. \rule{2mm}{2mm}

\section{Independence polynomials of well-covered spiders}

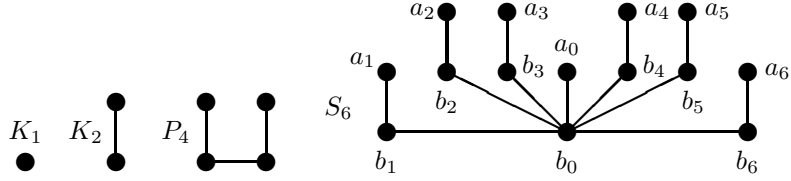
\begin{figure}[h]
\setlength{\unitlength}{0.8cm} 
\begin{picture}(5,3)\thicklines

  \put(2,0){\circle*{0.29}}
  \put(2,0.5){\makebox(0,0){$K_{1}$}}
  \put(3.5,0){\circle*{0.29}}
  \put(3.5,1){\circle*{0.29}} 
  \put(3.5,0){\line(0,1){1}}
  \put(3,0.5){\makebox(0,0){$K_{2}$}}  

  \multiput(5,0)(1,0){2}{\circle*{0.29}}
  \multiput(5,1)(1,0){2}{\circle*{0.29}}
  \put(5,0){\line(1,0){1}}
  \multiput(5,0)(1,0){2}{\line(0,1){1}}
  \put(4.5,0.5){\makebox(0,0){$P_{4}$}}

  \multiput(8,0.5)(3,0){3}{\circle*{0.29}}
  \multiput(8,1.5)(1,0){7}{\circle*{0.29}}
  \multiput(9,2.5)(1,0){2}{\circle*{0.29}}
  \multiput(12,2.5)(1,0){2}{\circle*{0.29}}
  \multiput(8,0.5)(1,0){6}{\line(1,0){1}}
  \multiput(8,0.5)(3,0){3}{\line(0,1){1}}
  \multiput(9,1.5)(1,0){2}{\line(0,1){1}}
  \multiput(12,1.5)(1,0){2}{\line(0,1){1}}
  \put(11,0.5){\line(-2,1){2}} 
  \put(11,0.5){\line(-1,1){1}}
  \put(11,0.5){\line(1,1){1}}
  \put(11,0.5){\line(2,1){2}}
  \put(11,0){\makebox(0,0){$b_{0}$}}
 \put(11,1.9){\makebox(0,0){$a_{0}$}}
 \put(8,0){\makebox(0,0){$b_{1}$}}
 \put(7.6,1.7){\makebox(0,0){$a_{1}$}}
 \put(9,1){\makebox(0,0){$b_{2}$}}
 \put(8.6,2.5){\makebox(0,0){$a_{2}$}}
 \put(10.45,1.6){\makebox(0,0){$b_{3}$}}
 \put(10.5,2.5){\makebox(0,0){$a_{3}$}}
 \put(12.45,1.6){\makebox(0,0){$b_{4}$}}
 \put(12.5,2.5){\makebox(0,0){$a_{4}$}}
 \put(13.1,1){\makebox(0,0){$b_{5}$}}
 \put(13.5,2.5){\makebox(0,0){$a_{5}$}}
 \put(14,0){\makebox(0,0){$b_{6}$}}
 \put(14.5,1.5){\makebox(0,0){$a_{6}$}}
 \put(7.2,0.9){\makebox(0,0){$S_{6}$}}

 \end{picture}
\caption{Well-covered spiders.}
\label{fig99}
\end{figure}

The well-covered spider $S_{n},n\geq 2$ has one vertex of degree $n+1,n$
vertices of degree $2$, and $n+1$ vertices of degree $1$ (see Figure \ref
{fig99}).

\begin{theorem}
The independence polynomial of any well-covered spider is unimodal,
moreover, 
\[
I(S_{n};x)=(1+x)\cdot \left\{ 1+\sum\limits_{k=1}^{n}\left[ {n \choose k}%
\cdot 2^{k}+ {n-1 \choose k-1}\right] \cdot x^{k}\right\} ,n\geq 2 
\]
and its mode is unique and equals $1+\left( n-1\right) \bmod 3+2\left(
\left\lceil n/3\right\rceil -1\right) $.
\end{theorem}

\setlength {\parindent}{0.0cm}\textbf{Proof.} Well-covered spiders comprise $%
K_{1},K_{2},P_{4}$ and $S_{n},n\geq 2$. Clearly, the independence
polynomials of $K_{1},K_{2},P_{4}$ are unimodal.%
\setlength
{\parindent}{3.45ex}

Using Proposition \ref{prop1}(i), we obtain the following formula for $S_{n}$%
: 
\begin{eqnarray*}
I(S_{n};x) &=&I(S_{n}-b_{0};x)+x\cdot I(S_{n}-N[b_{0}];x) \\
&=&(1+x)\cdot (1+2x)^{n}+x\cdot (1+x)^{n}=(1+x)\cdot R_{n}(x),
\end{eqnarray*}
where $R_{n}(x)=(1+2x)^{n}+x\cdot (1+x)^{n-1}$. By Lemma \ref{lem3}, to
prove that $I(S_{n};x)$ is unimodal, it is sufficient to show that $R_{n}(x)$
is unimodal. It is easy to see that 
\begin{eqnarray*}
R_{n}(x) &=&1+\sum\limits_{k=1}^{n}\left[ {n \choose k} \cdot 2^{k}+
{n-1 \choose k-1}\right] \cdot x^{k}=1+\sum\limits_{k=1}^{n}A(k)\cdot x^{k}, \\
A(k) &=& {n \choose k}\cdot 2^{k}+ {n-1 \choose k-1},1\leq k\leq n,A(0)=1.
\end{eqnarray*}

To start with, we show that $R_{n}$ is unimodal with the mode $k=n-1-\lfloor
(n-2)/3\rfloor $. Taking into account the proof of Lemma \ref{lem3}, namely,
the equality \ref{ModLe}, the mode of the polynomial 
\[
I(S_{n};x)=(1+x)\cdot R_{n}(x)=1+\sum\limits_{k=1}^{n}c_{k}\cdot x^{k} 
\]
is the index $m$, with $c_{m}=\max \{c_{k},c_{k+1}\}=\max
\{A(k)+A(k-1),A(k+1)+A(k)\}$. In other words, we will\emph{\ }give evidence
for 
\begin{equation}
m=\left\{ 
\begin{array}{l}
k,\ \ \ \ \ \ if \ A(k-1)>A(k+1), \\ 
k+1,\ if \ A(k-1)\leq A(k+1).
\end{array}
\right.  \label{mode}
\end{equation}

\begin{itemize}
\item  \textit{Claim 1.} If $n=3m+1$, then $R_{n}$ is unimodal with the mode 
$2m+1,I(S_{n};x)$ is also unimodal, and its unique mode equals $2m+1$.
\end{itemize}

We show that 
\begin{eqnarray*}
A(2m+i+1) &\geq &A(2m+i+2),0\leq i\leq m-1,\ and \\
A(2m+1-j) &\geq &A(2m-j),0\leq j\leq 2m.
\end{eqnarray*}

We have successively (for $2m+i+1=h$): 
\[
A(h)-A(h+1)=\left[ {3m+1 \choose h}\cdot 2^{h}+{3m \choose h-1}\right] -\left[ 
{3m+1 \choose h+1}\cdot 2^{h+1}+{3m \choose h}\right]
\]

\[
=\frac{(3m+1)!\cdot 2^{h}\cdot \left[ (h+1)-2\cdot (m-i)\right] }{%
(h+1)!\cdot (m-i)!}+\frac{(3m)!\cdot \left[ h-(m-i)\right] }{h!\cdot (m-i)!} 
\]

\[
=\frac{(3m+1)!\cdot \left( 3i+2\right) \cdot 2^{h}}{(h+1)!\cdot (m-i)!}+%
\frac{(3m)!\cdot \left( m+2i+1\right) }{h!\cdot (m-i)!}\geq 0. 
\]

Further, we get (for $2m-j=h$): 
\[
A(h+1)-A(h)=\left[ {3m+1 \choose h+1}\cdot 2^{h+1}+{3m \choose h}\right]
-\left[ {3m+1 \choose h}\cdot 2^{h}+ {3m \choose h-1}\right]
\]

\[
=\frac{(3m+1)!\cdot 2^{h}\cdot \left[ 2\cdot (m+j+1)-(h+1)\right] }{%
(h+1)!\cdot (m+j+1)!}+\frac{(3m)!\cdot \left[ (m+j+1)-h\right] }{h!\cdot
(m+j+1)!} 
\]

\[
=\frac{(3m)!}{(h+1)!\cdot (m+j+1)!}\left[ (3m+1)\cdot (3j+1)\cdot
2^{h}-(m-2j-1)\cdot (h+1)\right] \geq 0, 
\]
because $3m+1>m\geq m-2j-1$ and $2^{h}\geq h+1$.

To find the location of the mode of $I(S_{n};x)$, we obtain 
\[
A(2m)-A(2m+2)= 
\]

\[
={3m+1 \choose 2m}\cdot 2^{2m}+{3m \choose 2m-1}-{3m+1 \choose 2m+2}\cdot
2^{2m+2}+{3m \choose 2m+1} 
\]

\[
=\frac{(3m)!\cdot 2^{2m}}{(2m)!\cdot (m-1)!}\cdot \left[ \frac{3m+1}{m\cdot
(m+1)}\cdot (\frac{1}{2m}-\frac{1}{2m+1}\right] + 
\]

\[
+\,\frac{(3m)!}{(2m-1)!\cdot (m-1)!}\cdot \left[ \frac{1}{m\cdot (m+1)}-%
\frac{1}{2m\cdot (2m+1)}\right] >0. 
\]
Consequently, in accordance with equality (\ref{mode}), the mode of $%
I(S_{n};x)$ equals $2m+1$. Since $A(2m)-A(2m+2)>0$, the mode is unique.

\begin{itemize}
\item  \textit{Claim 2.} If $n=3m$, then $R_{n}$ is unimodal with the mode $%
2m$, $I(S_{n};x)$ is also unimodal,\emph{\ }and its unique mode equals $2m+1$%
.
\end{itemize}

We show that 
\begin{eqnarray*}
A(2m+i) &\geq &A(2m+i+1),0\leq i\leq m-1, \ and\ \\
A(2m-j) &\geq &A(2m-j-1),0\leq j\leq 2m-1.
\end{eqnarray*}

We have successively (for $2m+i=h$): 
\[
A(h)-A(h+1)=\left[ {3m \choose h}2^{h}+{3m-1 \choose h-1}\right] -\left[ 
{3m \choose h+1}2^{h+1}+{3m-1 \choose h}\right] 
\]

\[
=\frac{(3m)!\cdot 2^{h}\cdot \left[ (h+1)-2\cdot (m-i)\right] }{(h+1)!\cdot
(m-i)!}+\frac{(3m-1)!\cdot \left[ h-(m-i)\right] }{h!\cdot (m-i)!} 
\]

\[
=\frac{(3m)!\cdot (3i+1)\cdot 2^{h}}{(h+1)!\cdot (m-i)!}\cdot +\frac{%
(3m-1)!\cdot (m+2i)}{h!\cdot (m-i)!}\geq 0. 
\]

Further, we get (for $2m-j=h$): 
\[
A(h)-A(h-1)=\left[ {3m \choose h}\cdot 2^{h}+{3m-1 \choose h-1}\right] -\left[ 
{3m \choose h-1}\cdot 2^{h-1}+{3m-1 \choose h-2}\right] 
\]

\[
=\frac{(3m)!\cdot 2^{h-1}\cdot \left[ 2(m+j+1)-h\right] }{h!\cdot (m+j+1)!}+%
\frac{(3m-1)!\cdot \left[ (m+j+1)-(h-1)\right] }{(h-1)!\cdot (m+j+1)!} 
\]

\[
=\frac{(3m-1)!}{h!\cdot (m+j+1)!}\cdot \left[ 3m\cdot \frac{3j+2}{2}\cdot
2^{h}-(m-2j-2)\cdot h\right] \geq 0, 
\]
since $3m>m\geq m-2j-2$ and $2^{h}\geq h$.

To determine the mode of $I(S_{n};x)$, we obtain 
\[
A(2m-1)-A(2m+1)= 
\]

\[
={3m \choose 2m-1}\cdot 2^{2m-1}+{3m-1 \choose 2m-2}-{3m \choose 2m+1}\cdot
2^{2m+1}+{3m-1 \choose 2m} 
\]

\[
=\frac{3}{2}\cdot \frac{(3m-1)!}{(2m+1)!\cdot (m+1)!}\cdot \left[ (m-1)\cdot
(2m+1)-m\cdot 2^{2m}\right] <0. 
\]
Consequently, in accordance with equality (\ref{mode}), the mode of $%
I(S_{n};x)$ equals $2m+1$. Since $A(2m-1)-A(2m+1)<0$, the mode is unique.

\begin{itemize}
\item  \textit{Claim 3.} If $n=3m-1$, then $R_{n}$ is unimodal with the mode 
$2m-1$, $I(S_{n};x)$ is also unimodal, and its unique mode equals $2m$.
\end{itemize}

We show that 
\begin{eqnarray*}
A(2m+i-1) &\geq &A(2m+i),0\leq i\leq m-1,\ and \\
A(2m-j-1) &\geq &A(2m-j-2),0\leq j\leq 2m-2.
\end{eqnarray*}

We have successively (for $2m+i=h$): 
\[
A(h-1)-A(h)=\left[ {3m-1 \choose h-1}\cdot 2^{h-1}+{3m-2 \choose h-2}\right]
-\left[ {3m-1 \choose h}\cdot 2^{h}+{3m-2 \choose h-1}\right] 
\]

\[
=\frac{(3m-1)!\cdot 2^{h-1}\cdot \left[ h-2\cdot (m-i)\right] }{h!\cdot
(m-i)!}+\frac{(3m-2)!\cdot \left[ (h-1)-(m-i)\right] }{(h-1)!\cdot (m-i)!} 
\]

\[
=\frac{(3m-1)!\cdot 3i\cdot 2^{h-1}}{h!\cdot (m-i)!}+\frac{(3m-2)!\cdot
(m+2i-1)}{(h-1)!\cdot (m-i)!}\geq 0. 
\]

Further, we get (for $2m-j-1=h$): 
\[
A(h)-A(h-1)=\left[ {3m-1 \choose h}\cdot 2^{h}+{3m-2 \choose h-1}\right]
-\left[ {3m-1 \choose h-1}\cdot 2^{h-1}+{3m-2 \choose h-2}\right] 
\]

\[
=\frac{(3m-1)!\cdot 2^{h-1}\cdot \left[ 2\cdot (3m-h)-h\right] }{h!\cdot
(3m-h)!}+\frac{(3m-2)!\cdot \left[ (3m-h)-(h-1)\right] }{(h-1)!\cdot (3m-h)!}
\]

\[
=\frac{(3m-2)!}{h!\cdot (m+j+1)!}\cdot \left[ (3m-1)\cdot \frac{3j+3}{2}%
\cdot 2^{h}-h\cdot (m-2j-3)\right] \geq 0, 
\]
because $3m-1>m\geq m-2j-1$ and $2^{h}\geq h$.

Finally, we obtain that 
\[
A(2m-2)-A(2m)= 
\]

\[
={3m \choose 2m-2}\cdot 2^{2m-2}+{3m-1 \choose 2m-3}-{3m-1 \choose 2m}\cdot
2^{2m}+{3m-2 \choose 2m-1} 
\]

\[
=\frac{(3m-1)!}{(2m-1)!\cdot (m+1)!}\cdot \left[ m-2-3\cdot 2^{2m-2}\right]
<0. 
\]
Consequently, in accordance with equality (\ref{mode}), the mode of $%
I(S_{n};x)$ equals $2m$. Since $A(2m-1)-A(2m+1)<0$, the mode is unique. \rule%
{2mm}{2mm}

\section{Transformations of some well-covered trees}

If both $v_{1}$ and $v_{2}$ are vertices of degree at least two in $%
G_{1},G_{2}$, respectively, then $(G_{1};v_{1})\ominus (G_{2},v_{2})$
is an \textit{internal edge-join} of $G_{1},G_{2}$. Notice that the
edge-join of two trees is a tree, and also that two trees can be internal
edge-joined provided each one is of order at least three. An alternative
characterization of well-covered trees is the following:

\begin{theorem}
\cite{LevitMan1} A tree $T$ is well-covered if and only if either $T$ is a
well-covered spider, or $T$ is the internal edge-join of a number of
well-covered spiders.
\end{theorem}

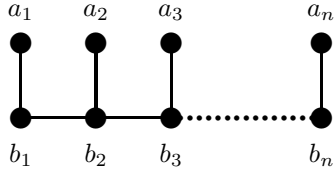
\begin{figure}[h]
\setlength{\unitlength}{1cm}%
\begin{picture}(5,2)\thicklines
  
  \multiput(4.5,0.5)(1,0){3}{\circle*{0.29}}
  \multiput(4.5,1.5)(1,0){3}{\circle*{0.29}}
  \put(8.5,0.5){\circle*{0.29}}
  \put(8.5,1.5){\circle*{0.29}}
  \put(4.5,0.5){\line(0,1){1}} 
  \put(5.5,0.5){\line(0,1){1}}
  \put(6.5,0.5){\line(0,1){1}}
  \put(4.5,0.5){\line(1,0){1}} 
  \put(5.5,0.5){\line(1,0){1}}
  \put(8.5,0.5){\line(0,1){1}}
  \multiput(6.5,0.5)(0.125,0){16}{\circle*{0.07}}
  \put(4.5,0){\makebox(0,0){$b_{1}$}}
  \put(4.5,1.9){\makebox(0,0){$a_{1}$}}
  \put(5.5,0){\makebox(0,0){$b_{2}$}}
  \put(5.5,1.9){\makebox(0,0){$a_{2}$}}
  \put(6.5,0){\makebox(0,0){$b_{3}$}}
  \put(6.5,1.9){\makebox(0,0){$a_{3}$}}
  \put(8.5,0){\makebox(0,0){$b_{n}$}}
  \put(8.5,1.9){\makebox(0,0){$a_{n}$}}
  
 \end{picture}
\caption{The centipede $W_{n}$.}
\label{fig1}
\end{figure}

As examples, $W_{n},n\geq 4$, and $G_{m,n},m\geq 2,n\geq 3$, (see Figures 
\ref{fig1} and \ref{fig999}) are internal edge-join of a number of
well-covered spiders, and consequently, they are well-covered trees. The aim
of this section is to show that the independence polynomials of $W_{n}$ and
of $G_{m,n}$ are unimodal. The idea is to construct, for these trees, some
claw-free graphs having the same independence polynomial, and then to use
Theorem \ref{th2}. We leave open the question whether the procedure we use
is helpful to define a claw-free graph with the same independence polynomial
as a general well-covered tree.

A \textit{centipede }is a tree\textit{\ }denoted by $W_{n}=(A,B,E),n\geq 1$,
(see\textit{\ }Figure \ref{fig1}), where $A\cup B$ is its vertex set, $%
A=\{a_{1},...,a_{n}\},B=\{b_{1},...,b_{n}\},A\cap B=\emptyset $, and the
edge set $E=\{a_{i}b_{i}:1\leq i\leq n\}\cup \{b_{i}b_{i+1}:1\leq i\leq
n-1\} $.

The following result was proved in \cite{LevitMan2}, but for the sake of
self-consistency of this paper and to illustrate the idea of a hidden
correspondence between well-covered trees and claw-free graphs, we give here
its proof.

\begin{theorem}
\cite{LevitMan2} For any $n\geq 1$ the following assertions hold:

(i) $I(W_{2n};x)=(1+x)^{n}Q_{n}(x)=I(\bigtriangleup _{n}\amalg nK_{1};x)$,
where $Q_{n}(x)=I(\bigtriangleup _{n};x)$;

$\quad I(W_{2n+1};x)=(1+x)^{n}Q_{n+1}(x)=I((K_{2}\ominus \bigtriangleup
_{n})\amalg nK_{1};x)$,

\quad where $Q_{n+1}(x)=I(K_{2}\ominus \bigtriangleup _{n};x)$;

(ii) $I(W_{n};x)$ is unimodal and 
\[
I(W_{n};x)=(1+x)\cdot (I(W_{n-1};x)+x\cdot I(W_{n-2};x)),n\geq 2, 
\]

\quad where $I(W_{0};x)=1,I(W_{1};x)=1+2x$.
\end{theorem}

\setlength {\parindent}{0.0cm}\textbf{Proof.} 
\setlength
{\parindent}{3.45ex} (i) Evidently, the polynomials $I(W_{n};x),1\leq n\leq
3 $, are unimodal, since 
\[
I(W_{1};x)=1+\mathbf{2}x,I(W_{2};x)=1+\mathbf{4}x+3x^{2},I(W_{3};x)=1+6x+%
\mathbf{10}x^{2}+5x^{3}. 
\]

Applying $\lfloor n/2\rfloor $ times Lemma \ref{lem2}(ii), we obtain that
for $n=2m\geq 4$, 
\[
I(W_{n};x)=I(\bigtriangleup _{m}\amalg mK_{1};x)=I(\bigtriangleup
_{m};x)\cdot (1+x)^{m}, 
\]
while for $n=2m+1\geq 5,$%
\[
I(W_{n};x)=I(K_{2}\ominus \bigtriangleup _{m}\amalg
mK_{1};x)=I(K_{2}\ominus \bigtriangleup _{m};x)\cdot (1+x)^{m}. 
\]

(iii) According to Proposition \ref{prop4} and Lemma \ref{lem3}, it follows
that $I(W_{n};x)$ is unimodal. Further, taking $U=\{a_{n},b_{n}\}$ and
applying Proposition \ref{prop1}(ii), we obtain: 
\begin{eqnarray*}
I(W_{n};x) &=&I(W_{n}-U;x)+x\cdot (I(W_{n}-N[a_{n}];x)+I(W_{n}-N[b_{n}];x))
\\
&=&I(W_{n-1};x)+x\cdot I(W_{n-1};x)+x\cdot (1+x)\cdot I(W_{n-2};x)) \\
&=&(1+x)\cdot (I(W_{n-1};x)+x\cdot I(W_{n-2};x)),
\end{eqnarray*}
which completes the proof. \rule{2mm}{2mm}\newline

It is worth mentioning that the problem of finding the mode of the centipede
is still unsolved.

\begin{conjecture}
\cite{LevitMan2} The mode of $I(W_{n};x)$ is $k=n-f(n)$ and $f(n)$ is given
by
\begin{eqnarray*}
f(n) &=&1+\lfloor n/5\rfloor ,2\leq n\leq 6, \\
f(n) &=&f(2+(n-2)\bmod5)+2\lfloor (n-2)/5\rfloor ,n\geq 7.
\end{eqnarray*}
\end{conjecture}

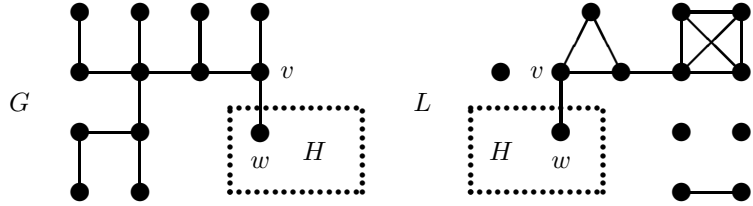
\begin{figure}[h]
\setlength{\unitlength}{0.8cm} 
\begin{picture}(5,4)\thicklines

  \multiput(2.5,0.5)(1,0){2}{\circle*{0.29}}
  \multiput(2.5,1.5)(1,0){2}{\circle*{0.29}}
  \multiput(2.5,2.5)(1,0){4}{\circle*{0.29}}
  \multiput(2.5,3.5)(1,0){4}{\circle*{0.29}}
  \multiput(2.5,0.5)(1,0){2}{\line(0,1){1}}
  \multiput(2.5,2.5)(1,0){3}{\line(0,1){1}}
  \put(3.5,1.5){\line(0,1){1}} 
  \put(2.5,1.5){\line(1,0){1}} 
  \put(2.5,2.5){\line(1,0){3}}
  \put(5.5,1.5){\circle*{0.29}}  
  \put(5.5,1.5){\line(0,1){2}} 
  \multiput(5,0.5)(0.15,0){15}{\circle*{0.07}}
 \multiput(5,1.9)(0.15,0){15}{\circle*{0.07}}
 \multiput(5,0.5)(0,0.15){10}{\circle*{0.07}}
 \multiput(7.2,0.5)(0,0.15){10}{\circle*{0.07}}
   \put(5.5,1){\makebox(0,0){$w$}}
   \put(6.4,1.2){\makebox(0,0){$H$}}
   \put(5.95,2.5){\makebox(0,0){$v$}}
   \put(1.5,2){\makebox(0,0){$G$}}

  \multiput(12.5,0.5)(1,0){2}{\circle*{0.29}}
  \multiput(9.5,2.5)(1,0){5}{\circle*{0.29}}
  \multiput(12.5,1.5)(1,0){2}{\circle*{0.29}}
  \multiput(12.5,3.5)(1,0){2}{\circle*{0.29}}
  \put(11,3.5){\circle*{0.29}}    
  \put(10.5,1.5){\circle*{0.29}} 
  \multiput(12.5,2.5)(1,0){2}{\line(0,1){1}}
  \put(12.5,2.5){\line(1,0){1}}
  \put(12.5,3.5){\line(1,0){1}}
  \put(12.5,2.5){\line(1,1){1}}
  \put(12.5,3.5){\line(1,-1){1}}
  \put(10.5,1.5){\line(0,1){1}} 
  \put(10.5,2.5){\line(1,0){3}} 
  \put(12.5,0.5){\line(1,0){1}}
  \put(10.5,2.5){\line(1,2){0.5}} 
  \put(11.5,2.5){\line(-1,2){0.5}}
\multiput(9,0.5)(0.15,0){15}{\circle*{0.07}}
 \multiput(9,1.9)(0.15,0){15}{\circle*{0.07}}
 \multiput(9,0.5)(0,0.15){10}{\circle*{0.07}}
 \multiput(11.2,0.5)(0,0.15){10}{\circle*{0.07}} 
   \put(10.5,1){\makebox(0,0){$w$}}
   \put(9.5,1.2){\makebox(0,0){$H$}}
   \put(10.1,2.5){\makebox(0,0){$v$}}
   \put(8.2,2){\makebox(0,0){$L$}}

 \end{picture}
\caption{The graphs $G$ and $L$ in Proposition 4.4(ii).}
\label{fig9}
\end{figure}

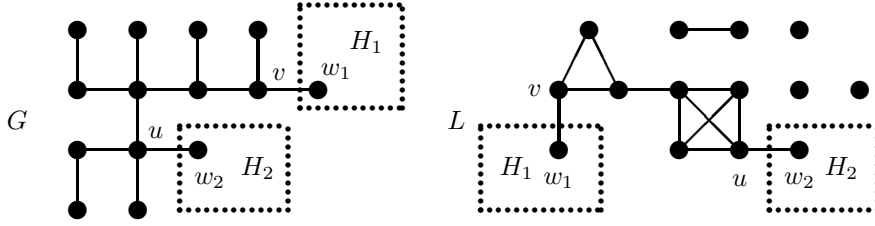
\begin{figure}[h]
\setlength{\unitlength}{0.8cm} 
\begin{picture}(5,4)\thicklines

  \multiput(1.5,0.5)(1,0){2}{\circle*{0.29}}
  \multiput(1.5,1.5)(1,0){2}{\circle*{0.29}}
  \multiput(1.5,2.5)(1,0){5}{\circle*{0.29}}
  \multiput(1.5,3.5)(1,0){4}{\circle*{0.29}}
  \multiput(1.5,0.5)(1,0){2}{\line(0,1){1}}
  \multiput(1.5,2.5)(1,0){4}{\line(0,1){1}}
  \put(2.5,1.5){\line(0,1){1}} 
  \put(1.5,1.5){\line(1,0){2}} 
  \put(1.5,2.5){\line(1,0){4}}
  \put(3.5,1.5){\circle*{0.29}}  
  \put(2.8,1.8){\makebox(0,0){$u$}}
 \multiput(3.2,0.5)(0.15,0){12}{\circle*{0.07}}
 \multiput(3.2,1.9)(0.15,0){12}{\circle*{0.07}}
 \multiput(3.2,0.5)(0,0.15){10}{\circle*{0.07}}
 \multiput(5,0.5)(0,0.15){10}{\circle*{0.07}}
  \put(3.7,1){\makebox(0,0){$w_{2}$}}
   
 \put(4.5,1.2){\makebox(0,0){$H_{2}$}}
 
  \put(4.85,2.75){\makebox(0,0){$v$}}
 
\multiput(5.2,2.2)(0.15,0){12}{\circle*{0.07}}
 \multiput(5.2,3.9)(0.15,0){12}{\circle*{0.07}}
 \multiput(5.2,2.2)(0,0.15){12}{\circle*{0.07}}
 \multiput(6.9,2.2)(0,0.15){12}{\circle*{0.07}}
 
  \put(5.8,2.8){\makebox(0,0){$w_{1}$}}
   \put(6.3,3.3){\makebox(0,0){$H_{1}$}}  
       \put(0.5,2){\makebox(0,0){$G$}}

  \multiput(11.5,3.5)(1,0){3}{\circle*{0.29}}
  \multiput(9.5,2.5)(1,0){6}{\circle*{0.29}}
  \multiput(11.5,1.5)(1,0){3}{\circle*{0.29}}
  
  \put(10,3.5){\circle*{0.29}}    
  \put(9.5,1.5){\circle*{0.29}} 
  \multiput(11.5,1.5)(1,0){2}{\line(0,1){1}}
  \put(9.5,1.5){\line(0,1){1}}
  \put(11.5,1.5){\line(1,0){2}}
  \put(9.5,2.5){\line(1,0){3}}
  \put(11.5,1.5){\line(1,1){1}}
  \put(11.5,2.5){\line(1,-1){1}}
  \put(11.5,3.5){\line(1,0){1}} 
\put(9.5,2.5){\line(1,2){0.5}}
  \put(10.5,2.5){\line(-1,2){0.5}} 
  
\multiput(8.2,0.5)(0.15,0){14}{\circle*{0.07}}
 \multiput(8.2,1.9)(0.15,0){14}{\circle*{0.07}}
 \multiput(8.2,0.5)(0,0.15){10}{\circle*{0.07}}
 \multiput(10.2,0.5)(0,0.15){10}{\circle*{0.07}} 
   \put(9.5,1){\makebox(0,0){$w_{1}$}}
   \put(8.8,1.2){\makebox(0,0){$H_{1}$}}
   \put(9.1,2.5){\makebox(0,0){$v$}}

    \put(13.5,1){\makebox(0,0){$w_{2}$}}
   \put(14.2,1.2){\makebox(0,0){$H_{2}$}}
   \put(12.5,1){\makebox(0,0){$u$}}
 \multiput(13,0.5)(0.15,0){12}{\circle*{0.07}}
 \multiput(13,1.9)(0.15,0){12}{\circle*{0.07}}
 \multiput(13,0.5)(0,0.15){10}{\circle*{0.07}}
 \multiput(14.8,0.5)(0,0.15){10}{\circle*{0.07}} 
      \put(7.8,2){\makebox(0,0){$L$}} 

 \end{picture}
\caption{The graphs $G$ and $L$ in Proposition 4.4(iii).}
\label{fig989}
\end{figure}

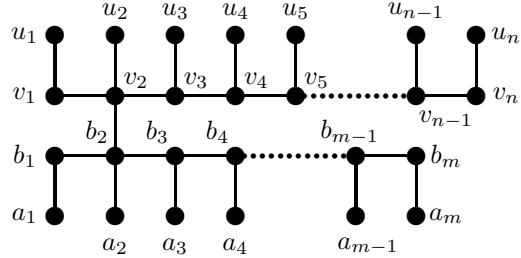
\begin{figure}[h]
\setlength{\unitlength}{0.8cm} 
\begin{picture}(5,4)\thicklines

  \multiput(4,0.5)(1,0){4}{\circle*{0.29}}
  \multiput(4,1.5)(1,0){4}{\circle*{0.29}}
  \multiput(9,0.5)(1,0){2}{\circle*{0.29}}
  \multiput(9,1.5)(1,0){2}{\circle*{0.29}}
  \multiput(4,2.5)(1,0){5}{\circle*{0.29}}
  \multiput(4,3.5)(1,0){5}{\circle*{0.29}}
  \multiput(10,2.5)(1,0){2}{\circle*{0.29}}
  \multiput(10,3.5)(1,0){2}{\circle*{0.29}}
 \multiput(7,1.5)(0.15,0){13}{\circle*{0.07}}
 \multiput(8,2.5)(0.15,0){13}{\circle*{0.07}}
  \multiput(4,0.5)(1,0){4}{\line(0,1){1}}
  \multiput(4,2.5)(1,0){5}{\line(0,1){1}}
  \multiput(9,0.5)(1,0){2}{\line(0,1){1}}
  \multiput(10,2.5)(1,0){2}{\line(0,1){1}}
  \put(5,1.5){\line(0,1){1}} 
  \put(4,1.5){\line(1,0){2}} 
  \put(4,2.5){\line(1,0){4}}
  \put(6,1.5){\line(1,0){1}} 
  \put(9,1.5){\line(1,0){1}} 
  \put(10,2.5){\line(1,0){1}}

  \put(3.5,1.5){\makebox(0,0){$b_{1}$}}
  \put(3.5,0.5){\makebox(0,0){$a_{1}$}}
  \put(4.7,1.9){\makebox(0,0){$b_{2}$}}
  \put(5,0){\makebox(0,0){$a_{2}$}}
  \put(5.7,1.9){\makebox(0,0){$b_{3}$}}
  \put(6,0){\makebox(0,0){$a_{3}$}}
  \put(6.7,1.9){\makebox(0,0){$b_{4}$}}
  \put(7,0){\makebox(0,0){$a_{4}$}}
  \put(8.9,1.9){\makebox(0,0){$b_{m-1}$}}
  \put(9.2,0){\makebox(0,0){$a_{m-1}$}}
  \put(10.5,1.5){\makebox(0,0){$b_{m}$}}
  \put(10.5,0.5){\makebox(0,0){$a_{m}$}}

  \put(3.5,2.5){\makebox(0,0){$v_{1}$}}
  \put(3.5,3.5){\makebox(0,0){$u_{1}$}}
  \put(5.35,2.75){\makebox(0,0){$v_{2}$}}
  \put(5,3.9){\makebox(0,0){$u_{2}$}}
  \put(6.35,2.75){\makebox(0,0){$v_{3}$}}
  \put(6,3.9){\makebox(0,0){$u_{3}$}}
  \put(7.35,2.75){\makebox(0,0){$v_{4}$}}
  \put(7,3.9){\makebox(0,0){$u_{4}$}}
  \put(8.35,2.75){\makebox(0,0){$v_{5}$}}
  \put(8,3.9){\makebox(0,0){$u_{5}$}}
  \put(10.5,2.1){\makebox(0,0){$v_{n-1}$}}
  \put(10,3.9){\makebox(0,0){$u_{n-1}$}}
  \put(11.5,2.5){\makebox(0,0){$v_{n}$}}
  \put(11.5,3.5){\makebox(0,0){$u_{n}$}}

 \end{picture}
\caption{The graph $G_{m,n},m\geq 2,n\geq 3$.}
\label{fig999}
\end{figure}

\begin{proposition}
\label{prop3}The following assertions are true:

(i) $I(G_{2,4};x)$ is unimodal, moreover $I(G_{2,4};x)=I(3K_{1}\amalg
K_{2}\amalg (K_{4}\ominus K_{3});x)$

\qquad (see Figure \ref{fig999});

(ii) $I(G;x)=I(L;x)$, where $G=(G_{2,4};v_{4})\ominus (H;w)$ and

$\qquad L=3K_{1}\amalg K_{2}\amalg (K_{4}\ominus K_{3})\ominus H$
(see Figure \ref{fig9});

\qquad if $w$ is simplicial in $H$, and $H$ is claw-free, then $I(G;x)$ is
unimodal;

(iii) $I(G;x)=I(L;x)$, where $G=(H_{1};w_{1})\ominus
(v;G_{2,4};u)\ominus (H_{2};w_{2})$ and

$\qquad L=3K_{1}\amalg K_{2}\amalg ((H_{1};w_{1})\ominus
(v;K_{3})\ominus (K_{4};u)\ominus (w_{2};H_{2})$ (see Figure \ref
{fig989});

\qquad if $w_{1},w_{2}$ are simplicial in $H_{1},H_{2}$, respectively, and $%
H_{1},H_{2}$ are claw-free,

\qquad then $I(G;x)$ is unimodal.
\end{proposition}

\setlength {\parindent}{0.0cm}\textbf{Proof.} (i)%
\setlength
{\parindent}{3.45ex} Using Proposition \ref{prop1}(iii) and the fact that $%
I(W_{4};x)=1+8x+21x^{2}+22x^{3}+8x^{4}=(1+x)^{2}(1+2x)(1+4x)$, we get that 
\begin{eqnarray*}
I(G_{2,4};x) &=&I(G_{2,4}-b_{2}v_{2};x)-x^{2}\cdot I(G_{2,4}-N\left(
b_{2}\right) \cup N\left( v_{2}\right) ;x) \\
&=&1+12x+55x^{2}+125x^{3}+\mathbf{150}x^{4}+91x^{5}+22x^{6},
\end{eqnarray*}
which is clearly unimodal. On the other hand, it is easy to check that 
\[
I(G_{2,4};x)=(1+x)^{3}(1+2x)(1+7x+11x^{2})=I(3K_{1}\amalg K_{2}\amalg
(K_{4}\ominus K_{3});x). 
\]

(ii) According to Proposition \ref{prop1}(iii), we infer that 
\begin{eqnarray*}
I(G;x) &=&I(G-vw;x)-x^{2}\cdot I(G-N(v)\cup N(w);x) \\
&=&I(G_{2,4};x)\cdot I(H;x)-x^{2}\cdot (1+x)\cdot I(H-N_{H}\left[ w\right]
;x)\cdot I(W_{4};x) \\
&=&I(G_{2,4};x)\cdot I(H;x)- \\
&&-x^{2}\cdot (1+x)^{3}\cdot (1+2x)\cdot (1+4x)\cdot I(H-N_{H}\left[
w\right] ;x).
\end{eqnarray*}

Figure \ref{fig9} shows the graphs $G$ and $L$. 

Let us denote $Q_{1}=$ $3K_{1}\amalg K_{2}\amalg K_{4},Q_{2}=3K_{1}\amalg
K_{2}\amalg (K_{4}\ominus K_{3}),L=(Q_{2};v)\ominus (H;w)$ and $%
e=vw$. Then, Proposition \ref{prop1}(iii) implies that 
\begin{eqnarray*}
I(L;x) &=&I(L-vw;x)-x^{2}\cdot I(L-N(v)\cup N(w);x) \\
&=&I(Q_{2};x)\cdot I(H;x)-x^{2}\cdot I(Q_{1};x)\cdot I(H-N_{H}\left[
w\right] ;x) \\
&=&I(G_{2,4};x)\cdot I(H;x)- \\
&&-x^{2}\cdot (1+x)^{3}\cdot (1+2x)\cdot (1+4x)\cdot I(H-N_{H}\left[
w\right] ;x).
\end{eqnarray*}
Consequently, $I(G;x)=I(L;x)$ holds.

In addition, if $w$ is simplicial in $H$, and $H$ is claw-free, then $L$ is
claw-free and, by Theorem \ref{th2}, $I(L;x)$ is unimodal. Hence, $I(G;x)$
is unimodal, as well.

(iii) Let $e=uw_{2}\in E(G)$. Then, according to Proposition \ref{prop1}%
(iii), we get that 
\begin{eqnarray*}
I(G;x) &=&I(G-uw_{2};x)-x^{2}\cdot I(G-N(u)\cup N(w_{2});x) \\
&=&I(H_{2};x)\cdot I((G_{2,4};v)\ominus (H_{1};w_{1});x)- \\
&&-x^{2}\cdot (1+x)^{2}\cdot (1+2x)\cdot I(H_{2}-N_{H_{2}}\left[
w_{2}\right] ;x)\cdot I((P_{4};v)\ominus (H_{1};w_{1});x).
\end{eqnarray*}
Now, Lemma \ref{lem2}(i) implies that 
\[
I(G;x)=I(H_{2};x)\cdot I((G_{2,4};v)\ominus (H_{1};w_{1});x)- 
\]
\[
-x^{2}\cdot (1+x)^{2}\cdot (1+2x)\cdot I(H_{2}-N_{H_{2}}\left[ w_{2}\right]
;x)\cdot I((K_{1}\sqcup K_{3};v)\ominus (H_{1};w_{1});x). 
\]
Figure \ref{fig989} shows the graphs $G$ and $L$. 

Let $e=uw_{2}\in E(L)$. Again by Proposition \ref{prop1}(iii), we infer that 
\begin{eqnarray*}
I(L;x) &=&I(L-uw_{2};x)-x^{2}\cdot I(L-N(u)\cup N(w_{2});x) \\
&&=I(H_{2};x)\cdot I((3K_{1}\amalg K_{2}\amalg (K_{4}\ominus
K_{3});v)\ominus (H_{1};w_{1});x)- \\
&&-x^{2}\cdot (1+x)^{3}\cdot (1+2x)\cdot I(H_{2}-N_{H_{2}}\left[
w_{2}\right] ;x)\cdot I((K_{3};v)\ominus (H_{1};w_{1});x).
\end{eqnarray*}

Further, Proposition \ref{prop3}(ii) helps us to deduce that 
\[
I(H_{2};x)\cdot I((G_{2,4};v)\ominus (H_{1};w_{1});x)=I(H_{2};x)\cdot
I((3K_{1}\amalg K_{2}\amalg (K_{4}\ominus K_{3});v)\ominus
(H_{1};w_{1});x). 
\]

Eventually, since 
\[
I((K_{1}\sqcup K_{3};v)\ominus (H_{1};w_{1});x)=(1+x)\cdot
I((K_{3};v)\ominus (H_{1};w_{1});x), 
\]
we obtain $I(G;x)=I(L;x)$.

If $w_{1},w_{2}$ are simplicial in $H_{1},H_{2}$, respectively, and $%
H_{1},H_{2}$ are claw-free, then $L$ is claw-free and, therefore, $I(L;x)$
is unimodal, by Theorem \ref{th2}. Hence, $I(G;x)$ is unimodal, too. \rule%
{2mm}{2mm}%

\begin{theorem}
The independence polynomial of $G_{m,n}=(W_{m};b_{2})\ominus
(W_{n};v_{2})$ is

unimodal, for any $m\geq 2,n\geq 2$.
\end{theorem}

\setlength {\parindent}{0.0cm}\textbf{Proof.} \textit{Case 1.} Suppose that $%
m=2,3,n=3,4$. The polynomial $I(G_{2,3};x)$ is unimodal, because 
\begin{eqnarray*}
I(G_{2,3};x) &=&(1+x)^{2}\cdot (1+2x)\cdot (1+6x+7x^{2}) \\
&=&1+10x+36x^{2}+\mathbf{60}x^{3}+47x^{4}+14x^{5}.
\end{eqnarray*}
According to Proposition \ref{prop3}(i), the independence polynomial of $%
G_{2,4}$ is unimodal and $I(G_{2,4};x)=I(3K_{1}\amalg K_{2}\amalg
(K_{4}\ominus K_{3});x)$.\setlength
{\parindent}{3.45ex}

Further, $I(G_{3,3};x)$ is unimodal, since 
\begin{eqnarray*}
I(G_{3,3};x) &=&I(G_{3,3}-v_{2}b_{2};x)-x^{2}\cdot I(G_{3,3}-N(v_{2})\cup
N(b_{2});x) \\
&=&I(W_{3};x)\cdot I(W_{3};x)-x^{2}\cdot (1+x)^{4} \\
&=&1+12x+55x^{2}+126x^{3}+\mathbf{154}x^{4}+96x^{5}+24x^{6}.
\end{eqnarray*}

Finally, $I(G_{3,4};x)$ is unimodal, because 
\begin{eqnarray*}
I(G_{3,4};x) &=&I(G_{3,4}-b_{1}b_{2};x)-x^{2}\cdot I(G_{3,4}-N(b_{1})\cup
N(b_{2});x) \\
&=&(1+2x)\cdot I(G_{2,4};x)-x^{2}\cdot (1+x)^{2}\cdot (1+2x)\cdot I(P_{4};x)
\\
&=&1+14x+78x^{2}+227x^{3}+\mathbf{376}x^{4}+357x^{5}+181x^{6}+38x^{7}.
\end{eqnarray*}

\begin{figure}[h]
\setlength{\unitlength}{1cm}%
\begin{picture}(5,2)\thicklines
  
  \multiput(1,0.5)(1,0){10}{\circle*{0.29}}
  \multiput(2,1.5)(1,0){4}{\circle*{0.29}}
  \put(6.5,1.5){\circle*{0.29}}
  \multiput(8,1.5)(1,0){3}{\circle*{0.29}}

  \put(12,0.5){\circle*{0.29}}
  \put(12,1.5){\circle*{0.29}}
  \put(4,0.5){\line(1,0){6}} 
  \multiput(3,0.5)(1,0){3}{\line(0,1){1}}  
  \multiput(8,0.5)(1,0){3}{\line(0,1){1}}
  \put(4,1.5){\line(1,0){1}}
  \put(4,0.5){\line(1,1){1}}
  \put(4,1.5){\line(1,-1){1}}
  \put(6,0.5){\line(1,2){0.5}} 
  \put(7,0.5){\line(-1,2){0.5}}   
  \put(12,0.5){\line(0,1){1}} 
  \multiput(10,0.5)(0.125,0){16}{\circle*{0.07}}
  
  \put(7,0){\makebox(0,0){$w$}}
  \put(8,0){\makebox(0,0){$v_{5}$}}
  \put(9,0){\makebox(0,0){$v_{6}$}}
  \put(10,0){\makebox(0,0){$v_{7}$}}
  \put(12,0){\makebox(0,0){$v_{n}$}}
  \put(8,1.9){\makebox(0,0){$u_{5}$}}
  \put(9,1.9){\makebox(0,0){$u_{6}$}}
  \put(10,1.9){\makebox(0,0){$u_{7}$}}
  \put(12,1.9){\makebox(0,0){$u_{n}$}}

 \end{picture}
\caption{The graph $L_{1}=3K_{1}\amalg K_{2}\amalg K_{4}\ominus
K_{3}\ominus W_{n-4}$.}
\label{fig1212}
\end{figure}
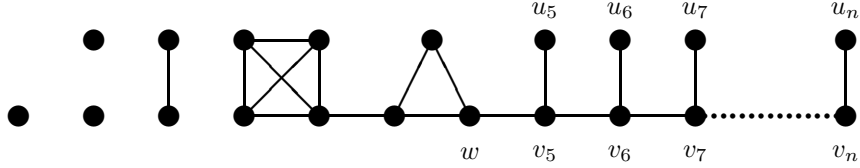

\textit{Case 2.} Assume that $m=2,n\geq 5$. According to Proposition \ref
{prop3} (ii), we infer that $I(G_{2,n};x)=I(L_{1};x)$, where $%
L_{1}=Q\ominus W_{n-4}$ and $Q=3K_{1}\amalg K_{2}\amalg
K_{4}\ominus K_{3}$ (see Figure \ref{fig1212}). 
Applying Lemma \ref{lem2}(ii), $I(L_{1};x)=I((mK_{1}\amalg (Q\ominus
mK_{3});x)$, if $n-4=2m$, and $I(L_{1};x)=I((mK_{1}\amalg (Q\ominus
mK_{3}\ominus K_{2});x)$, if $n-4=2m+1$. Since $mK_{1}\amalg
(Q\ominus mK_{3}\ominus K_{2})$ is claw-free, it follows that $%
I(L_{1};x)$ is unimodal, and consequently, $I(G_{2,n};x)$ is unimodal, too.

\textit{Case 3.} Assume that $m\geq 3,n\geq 5$. 
\begin{figure}[h]
\setlength{\unitlength}{1cm}%
\begin{picture}(5,3)\thicklines
  
 \multiput(1.5,1.5)(0.125,0){12}{\circle*{0.07}} 
  \multiput(1.5,0.5)(0,1){2}{\circle*{0.29}}
  \multiput(3,0.5)(1,0){8}{\circle*{0.29}}
  \multiput(3,1.5)(1,0){4}{\circle*{0.29}}
  \put(7.5,1.5){\circle*{0.29}}
  \multiput(9,1.5)(1,0){2}{\circle*{0.29}}
  \multiput(11.5,0.5)(0,1){2}{\circle*{0.29}}
  \put(5,0.5){\line(1,0){5}} 
  \put(3,1.5){\line(1,0){3}} 
  \put(1.5,0.5){\line(0,1){1}}
   \put(6.5,2.3){\line(1,0){1}}
   \multiput(6.5,2.3)(1,0){5}{\circle*{0.29}}
  \multiput(3,0.5)(1,0){4}{\line(0,1){1}}  
  \multiput(9,0.5)(1,0){2}{\line(0,1){1}}
  \put(5,1.5){\line(1,0){1}}
  \put(5,0.5){\line(1,1){1}}
  \put(5,1.5){\line(1,-1){1}}
  \put(7,0.5){\line(1,2){0.5}} 
  \put(8,0.5){\line(-1,2){0.5}}   
  \put(11.5,0.5){\line(0,1){1}} 
  \multiput(10,0.5)(0.125,0){12}{\circle*{0.07}}
  \put(4,1.9){\makebox(0,0){$b_{3}$}}
  \put(3,1.9){\makebox(0,0){$b_{4}$}}
  \put(1.5,1.9){\makebox(0,0){$b_{m}$}}
  \put(4,0){\makebox(0,0){$a_{3}$}}
  \put(3,0){\makebox(0,0){$a_{4}$}}
  \put(1.5,0){\makebox(0,0){$a_{m}$}}

  \put(9,0){\makebox(0,0){$v_{5}$}}
  \put(10,0){\makebox(0,0){$v_{6}$}}
  \put(11.6,0){\makebox(0,0){$v_{n}$}}

  \put(9.4,1.5){\makebox(0,0){$u_{5}$}}
  \put(10.4,1.5){\makebox(0,0){$u_{6}$}}
  \put(11.9,1.5){\makebox(0,0){$u_{n}$}}
  
 \end{picture}
\caption{The graph $L_{2}=W_{m-2}\ominus (3K_{1}\amalg K_{2}\amalg
K_{4}\ominus K_{3})\ominus W_{n-4}$.}
\label{fig121212}
\end{figure}
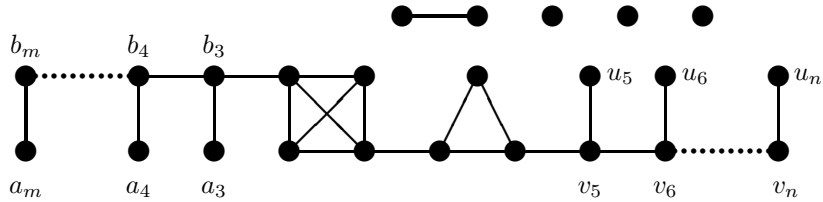
According to Proposition \ref{prop3}(iii), we obtain that $%
I(G_{m,n};x)=I(L_{2};x)$, where 
\[
L_{2}=W_{m-2}\ominus Q\ominus W_{n-4}\ and \ Q=3K_{1}\amalg
K_{2}\amalg K_{4}\ominus K_{3} 
\]
(see Figure \ref{fig121212}). Finally, by Theorem \ref{th2}, we infer that $%
I(L_{2};x)$ is unimodal, since by applying Lemma \ref{lem2}, $W_{m-2}$ and $%
W_{n-4}$ can be substituted by $pK_{1}\amalg (\ominus
pK_{3}\ominus K_{2})$ or $pK_{1}\amalg (\ominus pK_{3})$,
depending on the parity of the numbers $m-2,n-4$. Consequently, the
polynomial $I(G_{m,n};x)$ is unimodal, as well. \rule{2mm}{2mm}

\section{Conclusions}

In this paper we keep investigating the unimodality of independence
polynomials of some well-covered trees started in \cite{LevitMan2}. Any such
a tree is an edge-join of a number of ''\textit{atoms}'', called
well-covered spiders. We proved that the independence polynomial of any
well-covered spider is unimodal, straightforwardly indicating the location
of the mode. We also showed that the independence polynomial of some
edge-join of well-covered spiders is unimodal. In the later case, our
approach was indirect, via claw-free graphs.

\begin{figure}[h]
\setlength{\unitlength}{1cm}%
\begin{picture}(5,2)\thicklines
  
  \multiput(1,0)(1,0){2}{\circle*{0.29}}
  \multiput(1,1)(1,0){2}{\circle*{0.29}}
  \put(2,2){\circle*{0.29}}
  \put(1,0){\line(1,0){1}} 
  \put(1,1){\line(1,0){1}}
  \put(2,0){\line(0,1){2}}
  \put(0.5,1){\makebox(0,0){$H_{1}$}}

  \multiput(3.5,0)(0,1){3}{\circle*{0.29}}
  \multiput(4.5,0)(0,1){2}{\circle*{0.29}}
  \put(3.5,0){\line(1,0){1}} 
  \put(3.5,0){\line(0,1){1}}
  \put(3.5,1){\line(1,0){1}}
  \put(4.5,0){\line(0,1){1}}
  \put(3,1){\makebox(0,0){$H_{2}$}}

 \multiput(6,0)(0,1){3}{\circle*{0.29}}
  \multiput(7,0)(0,1){3}{\circle*{0.29}}
  \put(6,0){\line(1,0){1}} 
  \put(6,1){\line(1,0){1}}
  \put(6,2){\line(1,0){1}}
  \put(6,0){\line(0,1){2}}
  \put(7,0){\line(0,1){2}}
  \put(6,0){\line(1,1){1}}
  \put(6,1){\line(1,1){1}}
  \put(6,1){\line(1,-1){1}}
  \put(6,2){\line(1,-1){1}}
  \put(5.5,1){\makebox(0,0){$H_{3}$}}

  \put(10.5,0){\circle*{0.29}}
  \put(10.5,2){\circle*{0.29}}
  \multiput(9,1)(1,0){4}{\circle*{0.29}}
  \put(9,1){\line(1,0){3}} 
  \put(9,1){\line(3,2){1.5}}
  \put(10,1){\line(1,2){0.5}}
  \put(9,1){\line(3,-2){1.5}}
  \put(10,1){\line(1,-2){0.5}}
   \put(12,1){\line(-3,2){1.5}}
   \put(11,1){\line(-1,2){0.5}}
  \put(12,1){\line(-3,-2){1.5}}
  \put(11,1){\line(-1,-2){0.5}}
  \put(8.3,1){\makebox(0,0){$H_{4}$}}
  
  \end{picture}
\caption{$I(H_{1};x)=I(H_{2};x)$ and $I(H_{3};x)=I(H_{4};x)$.}
\label{fig70}
\end{figure}
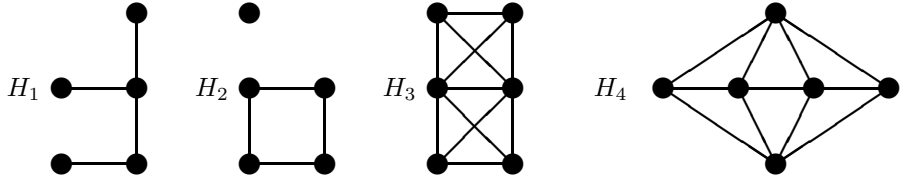

Let us notice that $I(H_{1};x)=I(H_{2};x)=1+5x+6x^{2}+2x^{3}$, and also $%
I(H_{3};x)=I(H_{4};x)=1+6x+4x^{2}$, where $H_{1},H_{2},H_{3},H_{4}$ are
depicted in Figure \ref{fig70}. In other words, there exist a well-covered
graph whose independence polynomial equals the independence polynomial of a
non-well-covered tree (e.g., $H_{2}$ and $H_{1}$), and also a well-covered
graph, different from a tree, namely $H_{4}$, satisfying $%
I(H_{3};x)=I(H_{4};x)$, where $H_{3}$ is not a well-covered graph. Moreover,
we can show that for any $\alpha \geq 2$ there are two connected graphs $%
G_{1},G_{2}$ such that $\alpha (G_{1})=\alpha (G_{2})=\alpha $ and $%
I(G_{1};x)=I(G_{2};x)$, but only one of them is well-covered. To see this,
let us consider the following two graphs: 
\[
G_{1}=L+(H_{1}\amalg H_{2}\amalg 2K_{1}),G_{2}=(L_{1}\amalg
L_{2})+(H_{1}\amalg H_{2}\amalg K_{2}), 
\]
where $L,L_{1},L_{2}$ are well-covered graphs, and 
\[
L=(L_{1},v_{1})\ominus
(L_{2},v_{2}),H_{1}=L_{1}-N[v_{1}],H_{2}=L_{2}-N[v_{2}],\alpha (L)=\alpha
(L_{1})+\alpha (L_{2}). 
\]
It follows that $\alpha (H_{1})=\alpha (L_{1})-1,\alpha (H_{2})=\alpha
(L_{2})-1$, and therefore, we obtain $\alpha (G_{1})=\alpha (G_{2})=\alpha
(L)$. It is easy to check that $G_{1}$ is well-covered, while $G_{2}$ is not
well-covered. According to Proposition \ref{prop1}(iii), we infer that 
\begin{eqnarray*}
I(L;x) &=&I(L-v_{1}v_{2};x)-x^{2}\cdot I(L-N_{L}(v_{1})\cup N_{L}(v_{2});x)
\\
&=&I(L_{1};x)\cdot I(L_{2};x)-x^{2}\cdot I(H_{1};x)\cdot I(H_{2};x)
\end{eqnarray*}
which we can write as follows: 
\[
I(L;x)+(1+x)^{2}\cdot I(H_{1};x)\cdot I(H_{2};x)=I(L_{1};x)\cdot
I(L_{2};x)+(1+2x)\cdot I(H_{1};x)\cdot I(H_{2};x) 
\]
or, equivalently, as 
\[
I(L;x)+I(2K_{1};x)\cdot I(H_{1};x)\cdot I(H_{2};x)=I(L_{1};x)\cdot
I(L_{2};x)+I(K_{2};x)\cdot I(H_{1};x)\cdot I(H_{2};x). 
\]
In other words, we get: 
\[
I(L+(2K_{1}\amalg H_{1}\amalg H_{2});x)=I(L_{1}\amalg L_{2}+(K_{2}\amalg
H_{1}\amalg H_{2});x), 
\]
i.e., $I(G_{1};x)=I(G_{2};x)$.

However, in some of our findings we defined claw-free graphs that
simultaneously are well-covered and have the same independence polynomials
as the well-covered trees under investigation. These results give an
evidence for the following conjecture.

\begin{conjecture}
If $T$ is a well-covered tree and $I(T;x)=I(G;x)$, then $G$ is well-covered.
\end{conjecture}

\end{document}